\documentclass[3p,12pt]{elsarticle}

\newtheorem{proposition}{Proposition}
\newtheorem{theorem}{Theorem}
\newtheorem{corollary}{Corollary}

\newtheorem{remark}{Remark}
\newtheorem{definition}{Definition}
\newtheorem{example}{Example}
\usepackage{graphics}
\usepackage{amsmath}
\usepackage{amssymb}
\usepackage{epsfig}
\usepackage{epstopdf}
\usepackage{color}
\usepackage[%
breaklinks%
,colorlinks%
,linkcolor=red
,anchorcolor=blue%
,pagecolor=blue%
,citecolor=blue%
,bookmarks=ture%
]{hyperref}

\def\qed{\hfill  \framebox(5,5){}}
\def\deg{{\rm deg}}
\def\cP{{\mathcal P}}

\def\index{{\rm index}}
\def\egcd{{\rm \epsilon gcd}}
\def\eindex{{\rm \epsilon index}}

\def\numer{{\rm numer}}

\def\deg{{\rm deg}}

\def\lc{{\rm lc}}

\def\cP{{\mathcal P}}

\def\para{\vspace{2 mm}}

\def\resultant{{\rm Res}}

\def\index{{\rm index}}

\def\lc{{\rm lc}}

\def\coeff{{\rm coeff}}

 \def\res{{\rm Res}}
  \def\num{{\rm num}}
\def\cP{{\mathcal P}}

 \def\QRGCD{{\rm QRGCD}}

\begin{document}

\begin{frontmatter}

\title{Numerical Reparametrization of Rational Parametric \\Plane Curves}

\author[SPD]{Sonia P\'erez-D\'{\i}az}
\address[SPD]{Dpto. de F\'{\i}sica y Matem\'aticas,
 Universidad de Alcal\'a,
      E-28871 Madrid, Spain}
      \ead{sonia.perez@uah.es}

\author[SLY]{Li-Yong Shen\corref{cor}}
 \address[SLY]{School of Mathematical Sciences, University of
 CAS, Beijing, China}
 \ead{shenly@amss.ac.cn}

 \cortext[cor]{Corresponding author}

\begin{abstract}
In this paper, we present an algorithm for reparametrizing
 algebraic plane curves from a numerical point of view. That is, we deal with mathematical objects that are assumed to be given approximately.  More precisely,  given a tolerance $\epsilon>0$ and a rational parametrization $\cal P$ with perturbed float coefficients of a plane curve $\cal C$, we present an
algorithm that computes a parametrization $\cal Q$ of a new plane curve $\cal D$ such that ${\cal Q}$ is an {\it $\epsilon$--proper reparametrization} of $\cal D$. In addition, the error bound is carefully discussed and we present a formula that  measures the ``closeness"  between the input curve $\cal C$ and the output curve
$\cal D$.
\end{abstract}

\end{frontmatter}


\section{Introduction}

A rational parametrization $\cal P$ of an algebraic plane  curve   $\cal C$
establishes a rational co\-rres\-pon\-den\-ce  of $\cal C$
  with the affine or projective line. This correspondence is a birational
equivalence if $\cal P$  is proper i.e., if $\cal P$  traces the curve once. Otherwise, if $\cal P$  is not proper, to almost all points on $\cal C$, there corresponds more than one parameter value. L\"{u}roth's theorem shows  constructively that it is always possible to reparametrize an improperly parametrized curve such that it becomes properly parametrized. That is, to almost all points $p\in {\cal C}$ (except perhaps finitely many) there corresponds exactly one parameter value $t_0\in {\Bbb C}$ such that ${\cal P}(t_0)=p$. A proper reparametrization always reduces the degree of the rational functions defining the curve.

\para

The reparametrization problem, in particular when the variety is a curve or a
surface,  is specially interesting in some practical applications
in  computer aided
geometric design (C.A.G.D) where objects are often given and manipulated
parametrically. In addition, proper parametrizations play an
important role in many practical applications in C.A.G.D, such as in visualization (see \cite{HSW},
\cite{HL97}) or rational parametrization of offsets (see
\cite{ASS}). Also, they provide an implicitization approach based
on resultants (see \cite{CLO2} and \cite{Sen2}). Hence, the study of proper reparametrization has been concerned by some authors such as~\cite{chionh06,Perez-repara,diaz02,Sen2,shen06}, and several  efficient proper reparametrization
algorithms can be found
in~\cite{gao92,Perez-repara,Sed86}.

\para

The problem  of proper reparametrization for curves  has been widely discussed in symbolic consideration. More precisely, {\it given the  field of complex numbers $\Bbb C$, and a rational
parametrization ${\cal P}(t)\in {\Bbb C}(t)^2$ of an algebraic
plane curve $\cal C$ with exact
coefficients, one finds a rational proper parametrization
${\cal Q}(t)\in {\Bbb C}(t)^2$ of $\cal C$, and a rational
function $R(t)\in {\Bbb C}(t)\setminus{\Bbb C}$ such that ${\cal
P}(t)={\cal Q}(R(t))$.} Nevertheless, in many practical applications, for instance in the frame of
C.A.G.D, these approaches tend to be insufficient, since in
practice most of data objects are given or become approximate. As a consequence, there has been an increasing interest for the development of hybrid
symbolic-numerical algorithms, and approximate algorithms.

\para

Intuitively speaking,
one is given a tolerance
$\epsilon>0$, and an irreducible affine algebraic plane curve
$\cal C$ defined by a parametrization $\cal P$  with perturbed float coefficients   that is ``{\it
nearly improper}"   (i.e. improper within the tolerance $\epsilon$), and the problem consists in computing a rational
curve ${\mathcal C}$ defined by a  parametrization $\cal Q$, such that $\cal Q$ is proper and almost all points of the rational curve ${\mathcal D}$ are in the ``{\it
vicinity}" of $\mathcal C$. The notion of vicinity may be introduced
as the offset region limited by the external and internal offset
to $\mathcal C$ at distance $\epsilon$ (see Section 4 for more
details), and
therefore the problem consists in finding, if it is possible, a
rational curve  ${\mathcal D}$ properly parametrized and lying within the offset
region of $\mathcal C$.  For instance, let us suppose that  we are
given a tolerance $\epsilon=0.2$, and a curve $\mathcal C$ defined by the parametrization
{\small \[{\cal P}(t)=\left(\frac{1.999 t^2+3.999 t+2.005-0.003t^4+0.001t^3}{2.005+0.998t^4+4.002t^3+6.004t^2+3.997t}, \frac{0.001-0.998t^4-4.003t^3-5.996t^2-4.005t}{2.005+0.998t^4+4.002t^3+6.004t^2+3.997t}\right).\]}


\begin{figure}[h]
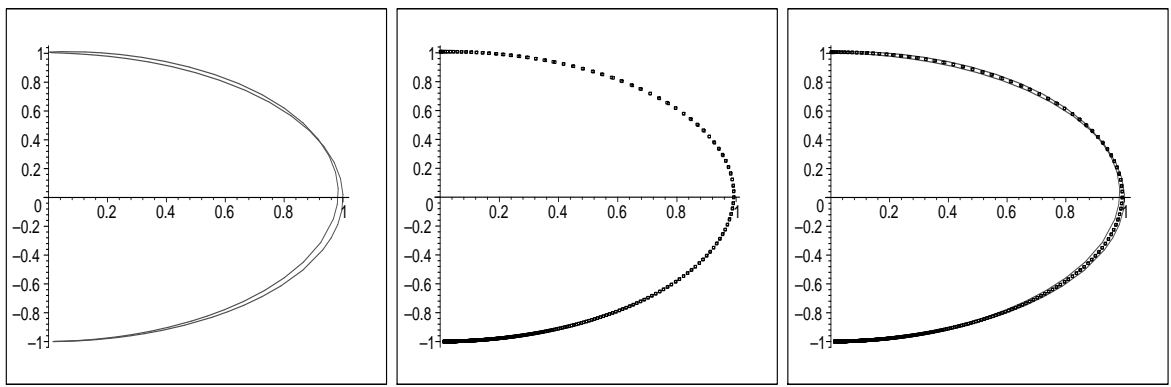

\begin{center}
\hbox{
\centerline{\epsfig{figure=example01.eps,width=5cm,height=5cm} \epsfig{figure=example02.eps,width=5cm,height=5cm} \epsfig{figure=example03.eps,width=5cm,height=5cm}
}}
\end{center}
\vspace*{-0.8cm}
\caption{Input curve $\mathcal C$ (left),  curve ${\mathcal D}$ (center),  curves $\mathcal C$ and ${\mathcal D}$ (right)}
 \label{Ex2fig1}
\end{figure}

One may check that $\cal P$ is proper in exact consideration
but it is nearly improper   since for almost all points $p:=\cP(s_0)\in {\cal C}$, there exist two values of the parameter $t$, given by the approximate roots   of the equation
$.4901606943 t^2+.2393271335\,10^{-8}(2202769 s_0+417838122)t-.4954325182 s_0^2-s_0=0,$
such that $\cP(t)$ is ``{\it almost equal}" to ${\cal P}(s_0)$.  Our method provides as
an answer  the curve ${\mathcal D}$ defined by the {\it ${\epsilon}$-proper reparametrization}
 $$  {\cal Q}(t)=  \left({\frac {- 0.00139214373770521\,{t}^{2}- 0.455587113115768\,t+
 0.230804565878748}{ 0.472790306463932\,{t}^{2}- 0.475516806696674\,t+
 0.233345983511073}}
,\,\right.$$$$ \left.{\frac {- 0.472791477433681\,{t}^{2}+ 0.473001908925789\,t-
 0.00421763512489261}{ 0.472790306463932\,{t}^{2}- 0.475516806696674\,
t+ 0.233345983511073}}
\right).$$
In Figure 1, one may check that $\mathcal C$ and $\cal D$ are ``{\it close}".

\para

The problem of relating the tolerance with the vicinity notion, may be approached either analyzing locally the condition number of the implicit equations (see \cite{Farouki}) or studying whether for almost every point $p$ on the original curve, there exists a point $q$ on the output curve such that the euclidean distance of $p$ and $q$ is significantly smaller than the tolerance. In this paper our error analysis will be based on the second approach. From this fact, and using \cite{Farouki}, one may derive upper bounds for the distance of the offset region.

\para

Approximate algorithms have been developed for some applied numerical topics, such as, approximate parametrization of algebraic curves and surfaces~\cite{PSS,PSS1,PSS3}, approximate greatest common divisor (gcd)~\cite{Beckermann1, Beckermann2, Corless, erich08, Karmarkar, zeng04}, finding zeros of multivariate systems~\cite{Corless}, and factoring polynomials~\cite{Cor2,galliao02}.  Few papers discussed the problem of properly
reparametrizing a given parametric curve with perturbed float coefficients.
As we know, only a heuristic algorithm was proposed  in \cite{Sed86}. However, the error analysis is not discussed, and no step is given to detect whether
a numerical curve is improper within a tolerance. In symbolic considerations, the tracing index is used to determine the properness of a parametrization of an algebraic plane curve (see \cite{Sen2, vdw72}). Essentially, it is the cardinality of a generic fibre of the parametrization, and from the geometric point of  view, the tracing index measures the number of times that a parametrization
traces a curve over the algebraic closure of the ground field. In this paper, we extend the concept to
the numerical situations that is, the approximate improper index is expected to be the number
of parameter value mapped in a neighborhood to a generic point of a given plane curve.
This gives the theoretical foundation for our further discussion.

\para

In this paper, we review  the symbolic algorithm of reparametrization for algebraic plane curves presented in \cite{Perez-repara}, and we generalize it for the numerical case. For this purpose, after we formally introduce the notion of approximate improper index, we define the equivalence of two numerical rational parametric curves.
The followed structure is similar to the symbolic situation, but the discussions are quite different. Some important properties are generalized to the numerical situation. Moreover, as the necessary work for the numerical discussion, the relation between the reparameterized and the original curve is subtly analyzed. As the error control, the approximate reparameterized curve obtained is restricted in the offset region of the original one (and reciprocally).

\para

More precisely, the paper is organized as follows. First, the symbolic algorithm  of proper reparameterization presented in \cite{Perez-repara} is briefly reviewed (see Section~2). In Section~3, the definition of approximate improper index ($\eindex$) and $\epsilon$-numerical reparametrization are proposed. In addition, we  construct the $\epsilon$-numerical reparametrization, and we prove that it is $\epsilon$-proper. Afterwards, we discuss the relation between the reparameterized curve and the input one, and we show the error analysis (see Section 4).
 In Section~5, the numerical algorithm is given as well as some examples. Finally, we conclude with
Section~6, where we propose topics for further study.


\section{Symbolic Algorithm of Reparametrization for Curves}

Before describing the method for the approximate case, and for reasons of completeness, in this section we briefly review some notions and the algorithmic approach to symbolically reparametrize  curves presented in \cite{Perez-repara}.
Let ${\Bbb C}$ be the  field of the complex numbers,
and $\cal C$ a rational algebraic plane curve
over $\Bbb C$. A parametrization $\cal P$ of $\cal C$ is proper if and only if the map
$${\cal P}:{\Bbb  C} \longrightarrow {\cal C} \subset {\mathbb C}^2; t\longmapsto {\cal
P}(t)$$
is birational, or equivalently, if for almost every point on $\cal C$ and for almost all values of the
parameter in $\Bbb C$ the mapping $\cal P$  is rationally bijective.
The notion of properness can also be stated algebraically in terms of fields of rational
functions. In fact, a rational parametrization $\cal P$  is proper if and only if the induced
monomorphism $\phi_{\cal P}$ on the fields of rational functions
$$ \phi_{\cal P}:{\mathbb C}({\cal C}) \longrightarrow  {\mathbb C}(t); R(x,y)\longmapsto R({\cal
P}(t))$$
is an isomorphism. Therefore, $\cal P$  is proper if and only if the mapping $\phi_{\cal P}$ is surjective,
that is, if and only if $\phi_{\cal P}({\Bbb C}({\cal C})) = {\mathbb C}({\cal P}(t)) = {\mathbb C}(t)$.
Thus, L\"{u}roth's theorem implies that
any rational curve over $\Bbb C$ can be properly parametrized  (see \cite{Abh},  \cite{Sen2}, \cite{van1}). Furthermore, given an improper parametrization, in
 \cite{Perez-repara}, \cite{Sed86} it is shown how to compute a new parametrization of the same curve
being proper.

\para

Intuitively speaking, we  say that  $\cal P$  is proper if and only if  ${\cal P}(t)$ traces $\cal C$ only once. In this sense, we may generalize the above notion by introducing the notion of tracing index of ${\cal P}(t)$. More precisely,  we say that $k\in {\Bbb N}$ is the tracing index of ${\cal P}(t)$, and we denote it by $\index({\cal P})$,
if all but finitely many points on $\cal C$ are generated, via ${\cal P}(t)$, by $k$ parameter values; i.e.
$\index({\cal P})$ represents the number of times that ${\cal P}(t)$ traces $\cal C$. Hence, the birationality of  $\phi_\cP$, i.e. the properness of $\cP(t)$, is
characterized  by tracing index 1 (for further details  see \cite{Sen2}).

 \para

 For reasons of completeness, we summarize some properties of the resultant that will be used throughout the paper. To start with,   we represent the  univariate resultant of two polynomials $A, B\in {\Bbb C}[x_1,\ldots,x_n,t]$  as $\res_t(A,B)$. The resultant over a commutative ring of two polynomials $A$ and $B$ is defined as  the determinant of the Sylvester matrix associated to $A$ and $B$. Thus, it holds that $\res_t(A,B)\in {\Bbb C}[x_1,\ldots,x_n]$, and $\res_t(A,B)=0$ if and only if $A, B$ have a common factor (depending on $t$). In addition,   the resultant is contained in the ideal generated by its two input
polynomials, and hence if $A(\alpha, b)=B(\alpha, b)=0$ where $\alpha=(\alpha_1,\ldots, \alpha_n)$, then $\res_t(A,B)(\alpha)=0$. Reciprocally, if  $\res_t(A,B)(\alpha)=0$, then $\lc(A,t)(\alpha)=\lc(B,t)(\alpha)=0$ ($\lc(A,t)$ denotes the leading coefficient of $A$ w.r.t.
$t$) or there exists $b\in {\Bbb C}$ such that $A(\alpha, b)=B(\alpha, b)=0$  (for more details see for instance Chapter 3 in \cite{Cox1998}, or Sections 5.8 and 5.9 in \cite{Vander}).

Additionally, we remind the reader the following specialization resultant property that will be used for our purposes:  if $\alpha\in {\Bbb C}^n$ is such that
$\deg_{t}(\varphi_{\alpha}(A))=\deg_{t}(A)$, and
$\deg_{t}(\varphi_{\alpha}(B))=\deg_{t}(B)-k$ then,
\[ \varphi_{\alpha}(\res_t(A,B))=\varphi_{\alpha}(\lc(A,t))^k \res_t(\varphi_{\alpha}(A),\varphi_{\alpha}(B)), \]
where $\varphi_{\alpha}$ is the  natural evaluation
homomorphism
$$\varphi_{\alpha}:   {\Bbb C}[x_1,\ldots,x_n,t]   \longrightarrow   {\Bbb C}[x_1,\ldots,x_n,t]; A(x_1,\ldots,x_n,t)   \longmapsto      A(\alpha_1,\ldots,\alpha_n,t)$$
 (see Lemma 4.3.1, pp.96 in \cite{win}).

 \para

 Finally,  given $R(t)=r_1(t)/r_2(t)\in {\Bbb C}(t)$, where $\gcd(r_1,\,r_2)=1$, we define   $\deg(R)$     as the maximum of $\deg(r_1)$ and $\deg(r_2)$.

\para

In the following, we outline the   algorithm developed in \cite{Perez-repara}  that computes a rational proper reparametrization of an improperly parametrized algebraic plane curve.   The algorithm is valid over any field, and it is based on the computation of polynomial gcds and univariate
resultants.

\noindent
\begin{center}
\fbox{\hspace*{2 mm}\parbox{16 cm}{ \vspace*{2 mm} {\bf Symbolic Algorithm
{\sf Reparametrization for Curves}.}

 \vspace*{2 mm}

\noindent
{\sc Input}: a rational affine parametrization ${\mathcal
P}(t)=\left({p_{1, 1}(t)}/{p_{1, 2}(t)},\,{p_{2, 1}(t)}/{p_{2,
2}(t)}\right)\in {\Bbb C}(t)^2,$
with $\gcd(p_{i,1}, p_{i,2})=1,\,i=1,2$,  of an algebraic plane curve $\cal C$.\\
{\sc Output}:  a  rational proper parametrization ${\cal Q}(t)\in {\Bbb C}(t)^2$ of
$\cal C$, and a rational function $R(t)\in {\Bbb C}(t)\setminus{\Bbb C}$ such that ${\cal
P}(t)={\cal Q}(R(t))$.
\begin{enumerate}
\item[1.] Compute  $H_{j}(t,s)=p_{j, 1}(t)p_{j, 2}(s)-p_{j, 1}(s)p_{j, 2}(t),\quad
j=1,2.$
\item[2.]  Determine the polynomial $S(t, s)=\gcd(H_1(t,s), H_2(t,s))=
C_m(t)s^m+\cdots+C_0(t).$
\item[3.]  If $\deg_t(S)=1$,  {\sc Return}   ${\cal Q}(t)={\mathcal P}(t)$, and
$R(t)=t$. Otherwise go to Step 4.
\item[4.]  Consider a rational function
$R(t) = \frac{C_{i}(t)}{C_{j}(t)} \in {\Bbb C}(t),$ such that
$C_{j}(t),\, C_{i}(t)$   are two of the polynomials obtained in Step 2 such that $\gcd(C_j, C_i)=1$,   and $C_jC_i\not\in {\Bbb C}$.
\item[5.]  For $k=1,2$, compute the polynomials
\[L_k(s, x_k)=\res_t(G_{k}(t, x_k), s C_{j}(t)-C_{i}(t))=(q_{k,2}(s)x_k-q_{k,1}(s))^{\deg(R)},\]
where $G_{k}(t, x_k)=x_kp_{k, 2}(t)-p_{k, 1}(t)$.
\item[6.] {\sc Return}   ${\cal
Q}(t)=\left({q_{1,1}(t)}/{q_{1,2}(t)},\,{q_{2,1}(t)}/{q_{2,2}(t)}\right) \in {\Bbb C}(t)^2,$
and $R(t)= {C_{i}(t)}/{C_{j}(t)}$.
\end{enumerate}}\hspace{2 mm}}
\end{center}

\para
\begin{remark}\label{R-eindex0}
It is proved that $\index({\cal P})=\deg_t(S)$ (see Theorem 2 in \cite{Sen2}). In addition, for all but finitely many values $\alpha$ of the variable $s$, $\deg_t(S)=\deg_t(\gcd(H_1(t,\alpha), H_2(t,\alpha)))$ (see Lemma 4 and Subsection 3.1 in \cite{Sen2}).
\end{remark}



\begin{example} Let $\cal C$ be the
rational curve defined by the parametrization
\[{\cal P}(t)=\left(\frac{p_{1,1}(t)}{p_{1,2}(t)},\,
\frac{p_{2,1}(t)}{p_{2,2}(t)}\right)=\]\[\left(\frac{10t^4+13t^3+17t^2+3t^5+24t+11+t^6}{(t^3+2)(t^2+3t+7+3t^3)},\, -\frac{2t^4-t^3-9t^2+5+t^6+t^5}{(t^3+2)^2}\right).\]
In Step 1 of the algorithm, we compute the polynomials\\

\noindent $H_{1}(t, s)=270t-270s+216t^2+39t^3+107t^4-19t^6+31t^5-26t^6s^3-49t^6s^2+26t^5s^3-31t^4s^2+91t^4s^3-11t^5s^2-195t^3s^2+195t^2s^3+234ts^3+19s^6-31s^5-107s^4-39s^3-216s^2-54ts^2+12ts^4+66ts^6+6ts^5+31t^2s^4+54t^2s+49t^2s^6+11t^2s^5-91t^3s^4-234t^3s+26t^3s^6-26t^3s^5-12t^4s+27t^4s^6+t^4s^5-27t^6s^4-66t^6s-8t^6s^5-t^5s^4-6t^5s+8t^5s^6,$\\

\noindent $H_{2}(t,s)=36t^2+24t^3-8t^4+t^6-4t^5-5t^6s^3-9t^6s^2-4t^5s^3-8t^4s^3-36t^3s^2+36t^2s^3-s^6+4s^5+8s^4-24s^3-36s^2+9t^2s^6+8t^3s^4+5t^3s^6+4t^3s^5-2t^4s^6+2t^6s^4+t^6s^5-t^5s^6.$\\

\noindent Now, we determine   $S(t, s)$. We obtain
$$S(t, s)=C_0(t)+C_1(t)s+C_2(t)s^2+C_3(t)s^3,$$
where $C_0(t)=-6t-2t^2+t^3$,\,\,\,$ C_1(t)=3t^3+6$,\,\,\,$ C_2(t)=t^3+2$,\, and $C_3(t)=-t^2-3t-1.$\\

\noindent Since $\deg_t(S)>1$, we go to Step 4 of
algorithm,  and we consider
\[R(t)= \frac{C_3(t)}{C_2(t)}=\frac{-t^2-3t-1}{t^3+2}.\]
Note  that $\gcd(C_2,\, C_3)=1$. Now, we compute the polynomials
\[L_1(s,x_1)=\res_t(G_{1}(t, x_1), sC_2(t)-C_3(t))=-961(-3x_1+sx_1+1-3s+s^2)^3,\]
\[L_2(s,x_2)=\res_t(G_{2}(t,x_2), sC_2(t)-C_3(t))=-961(-x_2-1+s+s^2)^3,\]
where $G_i(t, x_i)=x_ip_{i,2}(t)-p_{i,1}(t),\,i=1,2$ (see Step 5).
Finally, in Step 6, the algorithm outputs the proper
 parametrization ${\cal Q}(t)$, and the rational function $R(t)$
$${\cal Q}(t)=\left(-\frac{1-3t+t^2}{t-3},\,-1+t+t^2\right),\quad \quad\quad R(t)=\frac{-t^2-3t-1}{t^3+2}.$$
\end{example}


\section{The Problem of Numerical Reparametrization for Curves}

The problem  of numerical reparametrization for curves can be stated
as follows: {\sf given the field $\Bbb C$ of complex numbers, a tolerance $\epsilon>0$, and a rational
parametrization $${\cal
P}(t)=\left({p_{1,1}(t)}/{p_{1,2}(t)},\,{p_{2,1}(t)}/{p_{2,2}(t)}\right) \in {\Bbb C}(t)^2,$$ of an algebraic
plane curve $\cal C$ that is {\it approximate improper} (see Definition \ref{Def-approxindex}),  find a rational parametrization
${\cal Q}(t)\in {\Bbb C}(t)^2$ of an algebraic plane curve $\cal D$,
and a rational
function $R(t)\in {\Bbb C}(t)\setminus{\Bbb C}$ such that ${\cal Q}$ is an {\it $\epsilon$--proper reparametrization} of $\cal D$ (see Definition \ref{Def-repara}).  }

\para

In this section, the input and output are not assumed to be exact as in Section 2. Instead,  we  deal with mathematical objects that are
given approximately, probably because they proceed from an exact data that has been
perturbed under some previous measuring process or manipulation. Note that, in many practical applications, for instance in the frame of computer aided geometric
design,  most of data objects
are given or become approximate.


\para

In this new situation, the idea is to adapt the algorithm in Section 2 as follows. We consider a rational parametrization ${\cal P}(t)\in {\Bbb C}(t)^2$ of an algebraic
plane curve $\cal C$. We recall that because of a previous measuring process or manipulation, the parametrization ${\cal P}$ is assumed to be
given approximately.  Afterwards, one computes  the polynomials introduced in Step 1 of the symbolic algorithm presented in Section 2.

\para

In Step 2 of the symbolic algorithm, since we are working with mathematical objects that are assumed to be
given approximately, we have to compute the approximate $\gcd$, denoted by $\egcd$, instead of  the $\gcd$ (note that the gcd of two not exact input polynomials is always 1). There are different $\egcd$ algorithms proposed for inexact polynomials (see for instance, \cite{Beckermann1, Beckermann2, Corless, erich08, Karmarkar, zeng04}).  
Some typical algorithms of univariate polynomials are included in the mathematical softwares, for example, {\tt Maple} provides some $\egcd$ algorithms in the package SNAP.
We here introduce the $\egcd$ algorithm for a pair of univariate numeric polynomials by using QR factoring. It is implemented in {\tt Maple} as the function $\QRGCD$. The $\QRGCD(f, g, x, \epsilon)$ function returns univariate numeric polynomials $u, v ,d$ such that $d$ is an $\egcd$ for the input polynomials $f$ and $g$, and $u, v$ satisfy (with high probability)
$$\| u f + v g - d \|_2 < \epsilon\|(f,g,u,v,d)\|_2,\quad \| f - d f_1 \|_2 < \epsilon\|f\|_2,\quad \mbox{and}\quad \| g - d  g_1 \|_2 <\epsilon \|g\|_2,$$ where the polynomials $f_1$ and $f_2$ are cofactors of $f$ and $g$ with respect to the divisor $d$, and $\|(f,g,u,v,d)\|_2:=\max\{\|f\|_2, \|g\|_2,  \|u\|_2, \|v\|_2, \|d\|_2\}$.

\para

At this point, we need to generalize the concept of tracing index (see Section 2) to the numerical situations. For this purpose, in the following definition, we introduce the notion of {\it approximate improper index of $\cal P$}.

\begin{definition}\label{Def-approxindex}
We define  the
{\em approximate improper index} of $\cal P$ as  $\deg_t(S^{{\cal P}{\cal P}}_\epsilon),$ where $$S^{{\cal P}{\cal P}}_\epsilon(t,s)=\egcd(H^{{\cal P}{\cal P}}_1, H^{{\cal P}{\cal P}}_2),\quad H^{{\cal P}{\cal P}}_{j}(t,s)=p_{j, 1}(t)p_{j, 2}(s)-p_{j, 1}(s)p_{j, 2}(t),\,\,j=1,2,$$
and $s$ is a new variable. We denote it as $\eindex({\cal P})$. Furthermore, $\cal P$  is said to be {\em approximate improper} or {\em
$\epsilon$-improper} if $\eindex({\cal P})>1$. Otherwise, $\cal P$  is said to be {\em approximate  proper} or {\em
$\epsilon$-proper}.
\end{definition}


\para

Note that  in the symbolic situation, one can get the tracing index with probability one, by counting the common solutions  for a specialized $s_0$ (see Remark \ref{R-eindex0}).
For the numerical situation,  we can fix $s=s_0\in{\Bbb C}$ as a specialization and find the $\egcd$ for two univariate polynomials $H_1^{\cal P\cal P}(t,s_0)$ and $H_2^{\cal P\cal P}(t,s_0)$, under tolerance $\epsilon$.  Hence, we first can compute the approximate improper index by the specialization  and then, we can  recover an $\egcd$ defined by the polynomial  $S^{\cal P\cal P}_\epsilon(t,s)$. More precisely, $S^{{\cal P}{\cal P}}_\epsilon(t,s)$ can be found from several $S^{{\cal P}{\cal P}}_\epsilon(t,s_k),k=1,\ldots,n$, whose degrees equal to the approximate improper index. The polynomial $S^{{\cal P}{\cal P}}_\epsilon(t,s)$ can be computed using least squares method while $n$ is greater than the number of the unterminated coefficients (see the method presented in \cite{Sed86}). Note that  the approximate index is related to the selected $\epsilon$, and the used $\egcd$ algorithm.

\para

Once the polynomial $S^{{\cal P}{\cal P}}_\epsilon$ is computed, we consider the rational function $R(t)$ similarly as in Step 4 of the symbolic algorithm, and in Step 5 we compute the same resultant. Again, since we are working with approximate mathematical objects, the resultant does not factor as in the symbolic case. That is, if the input would have been an exact parametrization, the symbolic algorithm
  in Section 2 would  have output the parametrization  ${\cal
Q}(t)=\left({q_{1,1}(t)}/{q_{1,2}(t)},\,{q_{2,1}(t)}/{q_{2,2}(t)}\right) \in {\Bbb C}(t)^2,$ where
\[L_k(s, x_k)=\res_t(G_{k}(t, x_k), s C_{j}(t)-C_{i}(t))=(q_{k,2}(s)x_i-q_{k,1}(s))^{\deg(R)},\, \quad \mbox{and}\]
\[G_{k}(t, x_k)=x_kp_{k, 2}(t)-p_{k, 1}(t),\quad \mbox{for}\quad k=1,2.\]
However, in our case, ${q_{i,1}(s)}/{q_{i,2}(s)}$ are not exact roots of the polynomials $L_i(s, x_i)$ but
{\it $\epsilon$--roots} or {\it $\epsilon$--points} (see \cite{PSS}).  Thus, one may expect that a
small perturbation of $L_k$, provides a new polynomial that factorizes as above and the root of this new polynomial provides the
  parametrization.

\para

The notion of {\it  $\epsilon$-point} is introduced in \cite{PSS} as follows:  given a tolerance $\epsilon>0$, and a non-zero polynomial $A \in {\Bbb C}[t,s]$, we say that $(t_0, s_0)\in {\Bbb C}^2$ is an $\epsilon$-point of $A$, if  $|A(t_0,s_0))|\leq \epsilon \|A\|$, where $\|\cdot \|$ denotes the infinity norm, and  $| \cdot |$ is the absolute value (for further details in this notion see \cite{PSS}, \cite{PSS1}, \cite{PSS2}, \cite{PSS3}). We represent it as  $A(t_0,s_0)\approx_{\epsilon} 0$.  In Definition \ref{Def-operator}, we generalize this concept, and in particular  the operator ${\approx_{\epsilon}}$. For this purpose, in the following, $\num(\tau)$ represents the numerator of a rational function  $\tau(t)\in {\Bbb C}(t)$. 

\para

\begin{definition}\label{Def-operator}
Given two non-zero  polynomials $A_i \in {\Bbb C}[t,s]$ with $\|A_i \|=1,\,i=1,2$, we say that $A_1\approx_{\epsilon} A_2$, if  $\|A_1-A_2\|\leq \epsilon$ and $\deg_t(A_1)=\deg_t(A_2)$, $\deg_s(A_1)=\deg_s(A_2)$.   Furthermore, given $r(t)\in {\Bbb C}(t)$, and a non-zero polynomial  $A\in {\Bbb C}[t,s]$,  we say that $A(t,r(t))\approx_{\epsilon} 0$  if   $\|\num(A(t,r(t)))\|\leq \epsilon \|A\|$.
\end{definition}

 \para


\noindent
In the following, we consider a tolerance $\epsilon>0$, and $${\cal P}(t)=(p_1(t), p_2(t))=\left(\frac{p_{1, 1}(t)}{p_{1, 2}(t)}, \frac{p_{2, 1}(t)}{p_{2, 2}(t)}\right)\in
{\Bbb C}(t)^2,\quad \egcd(p_{j, 1}, p_{j, 2})=1,\,\,j=1,2$$  a rational parametrization of a given algebraic plane curve $\cal C$. We remind that ${\cal P}$ is expected to be given with perturbed float coefficients.  We assume   that $\index({\cal P})=1$. Observe that we are working numerically and then,   with
probability almost one $\deg_t(S)=1$, where $S$ is the polynomial introduced in Section 2. Otherwise, if $\index({\cal P})>1$, we may apply {\sf Symbolic Algorithm Reparametrization
for Curves} in Section 2.

\para

\noindent
We also consider the polynomials
$$S^{{\cal P}{\cal Q}}_\epsilon(t,s)=\egcd(H^{{\cal P}{\cal Q}}_1, H^{{\cal P}{\cal Q}}_2),\quad H^{{\cal P}{\cal Q}}_{j}(t,s)=p_{j, 1}(t)q_{j, 2}(s)-q_{j, 1}(s)p_{j, 2}(t),\,\,j=1,2,$$
where $s$ is a new variable, and
$${\cal Q}(t)=(q_1(t), q_2(t))=\left(\frac{q_{1, 1}(t)}{q_{1, 2}(t)}, \frac{q_{2, 1}(t)}{q_{2, 2}(t)}\right)\in
{\Bbb C}(t)^2,\quad \egcd(q_{j, 1}, q_{j, 2})=1,\,\,j=1,2$$   a rational parametrization of a new plane curve. Observe that these polynomials  generalize the polynomials introduced in Definition \ref{Def-approxindex}. In these conditions, we say that ${\cal P}(t) \sim_{\epsilon} {\cal Q}(r(t))$ if  $S^{{\cal P}{\cal Q}}_\epsilon(t,r(t)) \approx_{\epsilon} 0$,  where $r(t)\in {\Bbb C}(t)$ (see Definition \ref{Def-operator}). \\

 Observe that since $\egcd(p_{j, 1}, p_{j, 2})=1,\,\,j=1,2$, then $S^{{\cal P}{\cal Q}}_\epsilon(t,s)\in {\Bbb C}[t,s]\setminus{\Bbb C}[s]$. Similarly, since  $\egcd(q_{j, 1}, q_{j, 2})=1,\,\,j=1,2$, we also get that $S^{{\cal P}{\cal Q}}_\epsilon(t,s)\in {\Bbb C}[t,s]\setminus{\Bbb C}[t]$.\\

 \noindent
Throughout the paper, we  assume that   ${\cal P}(t) \not\sim_{\epsilon} (a,b)\in {\Bbb C}^2$ (see Remark \ref{R-eindex}, statement $1$).

\para

\begin{remark}\label{R-eindex} Observe that:
\begin{enumerate}
\item  Since ${\cal P}(t) \not\sim_{\epsilon} (a,b)\in {\Bbb C}^2$,  we have that $\deg_t(S^{{\cal P}{\cal P}}_\epsilon)\geq 1$. Indeed: note that $H^{{\cal P}{\cal P}}_j(t,s)\approx_{\epsilon} (t-s)N_j(t,s),$ where $N_j\in {\Bbb C}[t,s],\,\,j=1,2$. It holds that $N_j\not=0,\,j=1,2$; otherwise,  $S^{{\cal P}{\cal P}}_\epsilon(t,s)\approx_{\epsilon} 0$, and in particular $S^{{\cal P}{\cal P}}_\epsilon(t,s_0)\approx_{\epsilon} 0$ for $s_0\in {\Bbb C}$ satisfying that $p_{1,2}(s_0)p_{2,2}(s_0)\not=0$. Then, ${\cal P}(t) \sim_{\epsilon} {\cal P}(s_0)\in {\Bbb C}^2$ which is impossible, and thus $N_j\not=0,\,j=1,2$. Hence,   $S^{{\cal P}{\cal P}}_\epsilon(t,s)\approx_{\epsilon} (t-s)N(t,s)$, where  $N\in {\Bbb C}[t,s]\setminus \{0\}$.
\item   $\eindex({\cal P})=1$ if and only if $S^{{\cal P}{\cal P}}_\epsilon(t,s) \approx_\epsilon (t-s).$
\item Clearly the notion of {\em approximate improper index}  generalizes the notion of {\em  tracing index}. In particular, if $\eindex({\cal P})=1$ then $\index({\cal P})=1$.
\end{enumerate}
\end{remark}

\para

Now, we are ready to introduce the notions of {\it ${\epsilon}$-numerical reparametrization} and {\it ${\epsilon}$-proper reparametrization}.

\para

\begin{definition}\label{Def-repara}  Let ${\cal P}(t)=(p_1(t), p_2(t))\in
{\Bbb C}(t)^2$ be a rational parametrization of a given plane curve $\cal C$. We say that a parametrization ${\cal Q}(t)=(q_1(t), q_2(t))\in
{\Bbb C}(t)^2$ is an   ${\epsilon}$-numerical reparametrization  of ${\cal P}(t)$ if there exists $\displaystyle{{R}(t)={M(t)}/{N(t)}}\in {\Bbb C}(t)\setminus {\Bbb C},$  $\egcd(M, N)=1$, such that   ${\cal P} \sim_\epsilon {\cal Q}(R)$. In addition, if $\eindex({\cal Q})=1$, then we say that  ${\cal Q}$ is an    $\epsilon$-proper reparametrization  of ${\cal P}$. \end{definition}

\para

Using the concepts introduced above, we obtain some theorems where some properties characterizing numerical reparametrizations are proved. We start with Proposition \ref{Proposition-egcd}.

 \para

\begin{proposition}\label{Proposition-egcd} Let ${\cal Q}(t)=(q_1(t), q_2(t))\in
{\Bbb C}(t)^2,\,q_j=q_{j, 1}/q_{j, 2},\,\,\egcd(q_{j,1},q_{j,2})=1,\,\,j=1,2$ be such that  ${\cal Q}(t) \not\sim_{\epsilon} (a,b)\in {\Bbb C}^2$. Let   $\displaystyle{{R}(t)={M(t)}/{N(t)}}\in {\Bbb C}(t)\setminus {\Bbb C},$  $\egcd(M, N)=1$.  Up to constants in ${\Bbb C}\setminus\{0\}$, it holds that
\[S^{{\cal Q}(R){\cal Q}(R)}_\epsilon(t,s)\approx_\epsilon \num(S^{{\cal Q}{\cal Q}}_\epsilon(R(t),R(s))).\]
\end{proposition}
\noindent {\sc Proof.} From the definition of $S^{{\cal Q}{\cal Q}}_\epsilon(t,s)$, there are $M_1, M_2\in {\Bbb C}[t,s]$ satisfying that
\begin{equation}\label{eq1}
  H_j^{{\cal Q}{\cal Q}}(t,s)\approx_\epsilon  S^{{\cal Q}{\cal Q}}_\epsilon(t,s) M_j(t,s),\,\,j=1,2,\,\quad \mbox{and}\,\quad \egcd(M_1, M_2)=1
\end{equation}
(see Definition \ref{Def-operator}).
Now, taking into account the  definition of  $S^{{\cal Q}(R){\cal Q}(R)}_\epsilon$, one gets that
\begin{equation}\label{eq2}
S^{{\cal Q}(R){\cal Q}(R)}_\epsilon(t,s)=\egcd(H^{{\cal Q}(R){\cal Q}(R)}_1(t,s), H^{{\cal Q}(R){\cal Q}(R)}_2(t,s)).
\end{equation}
In addition, it holds that
{\small \begin{equation}\label{eq3}
\egcd(H^{{\cal Q}(R){\cal Q}(R)}_1(t,s), H^{{\cal Q}(R){\cal Q}(R)}_2(t,s))=\egcd(\num(H^{{\cal Q}{\cal Q}}_1(R(t),R(s))), \num(H^{{\cal Q}{\cal Q}}_2(R(t),R(s)))).
\end{equation}}
In order to prove \eqref{eq3}, we assume that   $\deg(q_{i,2})\geq \deg(q_{i,1})$ (otherwise, we reason similarly), and we consider  $q^*_{i,j}(x,y)\in {\Bbb C}[x,y]$  the homogenization of the polynomial  $q_{i,j}(x)\in {\Bbb C}[x]$ w.r.t. the variable $x$, and $\alpha:=\deg(q_{i,2})- \deg(q_{i,1})$. Under these conditions, equality \eqref{eq3} follows since
\[H^{{\cal Q}(R){\cal Q}(R)}_i(t,s)=\num\left(\frac{q_{i, 1}(R(t))}{q_{i, 2}(R(t))}-\frac{q_{i, 1}(R(s))}{q_{i, 2}(R(s))}\right)=\]\[N(t)^{\alpha}q^*_{i, 1}(M(t),N(t))q^*_{i, 2}(M(s),N(s))
 -N(s)^{\alpha}q^*_{i, 1}(M(s),N(s))q^*_{i, 2}(M(t),N(t)), \]
 and
 \[\num(H^{{\cal Q}{\cal Q}}_i(R(t),R(s)))=\num\left(q_{i, 1}(R(t))q_{i, 2}(R(s))-q_{i, 1}(R(s))q_{i, 2}(R(t)) \right)=\]
 \[\numer\left(\frac{q^*_{i, 1}(M(t),N(t))q^*_{i, 2}(M(s),N(s))}{N(t)^{\deg(q_{i,1})}N(s)^{\deg(q_{i,2})}}
 -\frac{q^*_{i, 1}(M(s),N(s))q^*_{i, 2}(M(t),N(t))}{N(t)^{\deg(q_{i,2})}N(s)^{\deg(q_{i,1})}}\right)=\]
 \[N(t)^{\alpha}q^*_{i, 1}(M(t),N(t))q^*_{i, 2}(M(s),N(s))
 -N(s)^{\alpha}q^*_{i, 1}(M(s),N(s))q^*_{i, 2}(M(t),N(t)). \]

 \noindent
Thus, from the above equalities, one deduces that
\[S^{{\cal Q}(R){\cal Q}(R)}_\epsilon(t,s)\overbrace{=}^{\mbox{\eqref{eq2} and \eqref{eq3}}}\egcd(\num(H^{{\cal Q}{\cal Q}}_1(R(t),R(s))), \num(H^{{\cal Q}{\cal Q}}_2(R(t),R(s))))\overbrace{\approx_\epsilon}^{\mbox{\eqref{eq1}}}\]\[ \egcd(\num(S^{{\cal Q}{\cal Q}}_\epsilon(R(t),R(s))) \num(M_1(R(t),R(s))), \num(S^{{\cal Q}{\cal Q}}_\epsilon(R(t),R(s))) \num(M_2(R(t),R(s))))\]\[ =\num(S^{{\cal Q}{\cal Q}}_\epsilon(R(t),R(s))) M(t,s),\]
where $M(t,s):=\egcd(\num(M_1(R(t),R(s))),\num(M_2(R(t),R(s)))).$ Since $\egcd(M_1, M_2)=1$, we have that $M(t,s)=1$, and we conclude that \[\hspace*{3cm}S^{{\cal Q}(R){\cal Q}(R)}_\epsilon(t,s)\approx_\epsilon \num(S^{{\cal Q}{\cal Q}}_\epsilon(R(t),R(s))).\hspace*{6cm}\mbox{\qed} \]

\para

\para

In the following, we consider
${\cal Q}(t)\in {\Bbb C}(t)^2$   an $\epsilon$-numerical reparametrization  of ${\cal P}(t)$, and     ${\cal P} \sim_\epsilon {\cal Q}(R)$ where $\displaystyle{{R}(t)}=M(t)/N(t)\in {\Bbb C}(t)\setminus {\Bbb C}$  (see Definition \ref{Def-repara}).
In these conditions, we have the following results.

 \para

 \begin{theorem}\label{C-egcd0}  ${\cal Q}$ is  $\epsilon$-proper if and only if $\eindex({\cal P})=\deg(R).$
 \end{theorem}
\noindent {\sc Proof.}
If $\eindex({\cal Q})=1$,  from  Proposition \ref{Proposition-egcd} and Remark \ref{R-eindex} (statement $2$),  one gets that
$$S^{{\cal P}{\cal P}}_\epsilon(t,s) \approx_\epsilon \num(S^{{\cal Q}{\cal Q}}_\epsilon(R(t),R(s))) \approx_\epsilon \num(R(t)-R(s))=M(t)N(s)-M(s)N(t).$$
Therefore,
$\eindex({\cal P})=\deg_t(S^{{\cal P}{\cal P}}_\epsilon)=\deg(R)$ (see Definitions \ref{Def-approxindex} and \ref{Def-operator}).
Reciprocally,   from Proposition \ref{Proposition-egcd}, we have that
   $$S^{{\cal P}{\cal P}}_\epsilon(t,s) \approx_\epsilon \num(S^{{\cal Q}{\cal Q}}_\epsilon(R(t),R(s)))$$  which implies that $\deg_t(S^{{\cal P}{\cal P}}_\epsilon)=\deg_t(S^{{\cal Q}{\cal Q}}_\epsilon)\deg(R)$ (see Definition \ref{Def-operator}). Therefore, if $\eindex({\cal P})=\deg(R)$, then $\deg_t(S^{{\cal Q}{\cal Q}}_\epsilon)=1$ and thus  ${\cal Q}$ is  $\epsilon$-proper (see Definition \ref{Def-approxindex}). \qed

\para

 \begin{corollary}\label{Th-eindex}   It holds that $\eindex({\cal P})=\eindex({\cal Q})\deg(R).$  \end{corollary}
\noindent {\sc Proof.} Reasoning as in proof  of Theorem \ref{C-egcd0},  one deduces that  $\deg_t(S^{{\cal P}{\cal P}}_\epsilon)=\deg_t(S^{{\cal Q}{\cal Q}}_\epsilon)\deg(R)$. Thus, from Definition \ref{Def-approxindex}, we conclude that  $\eindex({\cal P})=\eindex({\cal Q})\deg(R).$ \qed

\para

\begin{corollary}\label{Th-egcd} ${\cal Q}$ is  $\epsilon$-proper if and only if
$$S^{{\cal P}{\cal P}}_\epsilon(t,s)\approx_\epsilon \num(R(t)-R(s))=M(t)N(s)-M(s)N(t).$$  \end{corollary}
\noindent {\sc Proof.}
If $\eindex({\cal Q})=1$,  reasoning as in proof  of Theorem \ref{C-egcd0},  one deduces that  $S^{{\cal P}{\cal P}}_\epsilon(t,s)\approx_\epsilon M(t)N(s)-M(s)N(t).$
Reciprocally, if $S^{{\cal P}{\cal P}}_\epsilon(t,s)\approx_\epsilon \num(R(t)-R(s))$,    we get that $\eindex({\cal P})=\deg_t(S^{{\cal P}{\cal P}}_\epsilon)=\deg(R)$ (see Definition \ref{Def-operator}). Thus,  Corollary \ref{Th-eindex} implies that ${\cal Q}$ is  $\epsilon$-proper.  \qed

 \para

\begin{center}
\underline{\sf Construction of the Rational Function $R$}
\end{center}

\para

In the following, we construct a rational function $R(t)\in {\Bbb C}(t)\setminus {\Bbb C}$, such that there exists an $\epsilon$-proper reparametrization  of ${\cal P}$. That is, there exists ${\cal Q}$ such that ${\cal P}\sim_\epsilon {\cal Q}(R)$ and $\cal Q$ is $\epsilon$-proper (see Theorem \ref{Th-resultant} and Corollary \ref{C-resultant}). Hence, we are addressing the existence of the $\epsilon$-proper reparameterization.

\para

\noindent
For this purpose, we first write
\begin{equation}\label{constructR}
  S^{{\cal P}{\cal P}}_\epsilon(t,s)= C_m(t)  s^m+C_{m-1}(t)  s^{m-1}+\cdots+ C_{0}(t).
\end{equation}
This polynomial is computed from the input parametrization $\cal P$, and then it is known.  Furthermore,  taking into account Corollary \ref{Th-egcd}, we have that
\[S^{{\cal P}{\cal P}}_\epsilon(t,s) \approx_\epsilon \num(R(t)-R(s)),\]
where  $R(t)=M(t)/N(t)\in {\Bbb C}(t)\setminus {\Bbb C}$ is the unknown rational function we are looking for. That is, we look for $R$ satisfying the above condition.

\para

\noindent
In the symbolic situation,  Lemma 3 in \cite{Perez-repara} states that, up to constants in ${\Bbb C}\setminus\{0\}$,
 $$\num\left(\frac{C_{i}(t)}{C_{j}(t)}-\frac{C_{i}(s)}{C_{j}(s)}\right)=C_m(t)  s^m+C_{m-1}(t)  s^{m-1}+\cdots+ C_{0}(t), $$ where $C_i, C_j$ are such that $C_{i}C_{j}\not\in {\Bbb C}$, and   $\gcd(C_{i},\,C_{j})=1$.
\para

\noindent
Therefore,   the unknown rational function, $R(t)\in {\Bbb C}(t)\setminus {\Bbb C}$,  can  be constructed as
\[R(t)=\frac{C_{i}(t)}{C_{j}(t)}\in {\Bbb C}(t)\setminus {\Bbb C},\]
where $C_i$ and $C_j$ are from~\eqref{constructR} satisfying that:
\begin{itemize}
  \item[1)] $C_{i}C_{j}\not\in {\Bbb C}$,
  \item[2)] $\egcd(C_{i},\,C_{j})=1$,
  \item[3)] $S^{{\cal P}{\cal P}}_\epsilon(t,s) \approx_\epsilon \num(R(t)-R(s)).$
\end{itemize}

\para

\begin{center}
\underline{\sf Construction  of the $\epsilon$-Numerical Reparametrization ${\cal Q}$}
\end{center}

\para

In the following,  we consider the rational function $R(t)=\frac{C_{i}(t)}{C_{j}(t)}\in {\Bbb C}(t)\setminus {\Bbb C}$  computed as above. Then, we have that   $S^{{\cal P}{\cal P}}_\epsilon(t,s)\approx_\epsilon \num(R(t)-R(s))$. Hence, from Corollary  \ref{Th-egcd}, if   ${\cal Q}$ is  such that ${\cal P}\sim_\epsilon {\cal Q}(R)$, then  ${\cal Q}$ is   $\epsilon$-proper.

\para
\para

\noindent
In Theorem \ref{Th-resultant},  we show how to compute the $\epsilon$-numerical   reparametrization  ${\cal Q}$.


\para

\begin{theorem}\label{Th-resultant} Let
\[L_k(s, x_k)=\res_t(G_{k}(t, x_k), s C_{j}(t)-C_{i}(t)), \,\,\mbox{where}\,\,G_{k}(t, x_k)=x_kp_{k, 2}(t)-p_{k, 1}(t),\,\,k=1,2.\]
If  \[L_k(s, x_k)= (x_kq_{k,2}(s)-q_{k,1}(s))^{\ell}+\epsilon^{\ell} W_k(s, x_k),\,\,\,\,\,\|\num(W_k(R,p_k))\|\leq \|{H}^{\cP{\cal Q}}_{k}\|^\ell,\,\,\,k=1,2,\]
where $\ell:=\deg(R)$, and  $\egcd(q_{k,1}, q_{k,2})=1$,
then
 ${\cal
Q}(s)=\left(\frac{q_{1,1}(s)}{q_{1,2}(s)},\,\frac{q_{2,1}(s)}{q_{2,2}(s)}\right)$  is  an  $\epsilon$-numerical   reparametrization  of ${\cal P}$.
\end{theorem}
\noindent {\sc Proof.} First, we observe that $L_k\not=0$ (otherwise, $G_k$ and $s C_{j}(t)-C_{i}(t)$ have a common factor depending on $t$, which is impossible because $\gcd(C_i, C_j)=1$). In addition, it holds that
$\deg_{x_k}(L_k)=\deg(R)$. Indeed, since
\[L_k(s, x_k)=\res_t(G_{k}(t, x_k), s C_{j}(t)-C_{i}(t)),\]
we get that, up to constants in ${\Bbb C}(s)\setminus\{0\}$, $$L_k(s, x_k)=\prod_{\{\alpha_\ell\,|\,s C_{j}(\alpha_\ell)-C_{i}(\alpha_\ell)=0\}} G_k(\alpha_{\ell}, x_k),$$  (see  Sections 5.8 and 5.9 in \cite{Vander}), and thus $$\deg_{x_k}(L_k)=\deg_t(s C_{j}(t)-C_{i}(t))\deg_{x_k}(G_{k}(t, x_k))=\deg(R).$$ In addition, from $\deg_{x_k}(L_k)=\deg(R)$, we deduce that $\deg_{x_k}(W_k)\leq \ell$. In fact, since we are working numerically, we may assume w.l.o.g that $\deg_{x_k}(W_k)= \ell$.\\
Now, taking into account the properties of the resultant (see Section 2), one has that
\[0= L_k(R(t),p_k(t))= (p_k(t)q_{k,2}(R(t))-q_{k,1}(R(t)))^{\ell}+\epsilon^\ell W_k(R(t),p_k(t)).\] Then,   $$\num({H}^{\cP{\cal Q}}_{k}(t,R(t)))^\ell=\epsilon^\ell e_k(t),\,\,\, \mbox{where}\,\,\,e_k:=-W_k(R(t),p_k(t))p_{k,2}(t)^{\ell}C_j^{\ell\,\deg(q_k)},\qquad k=1,2.$$
Since $\deg_{x_k}(W_k)=\ell$, and  $\deg_{s}(W_k) = \ell\,\deg(q_k)$ (see Corollary \ref{C-degree}),  one has that  $$e_k=-\num(W_k(R(t),p_k(t)))\in {\Bbb C}[t]$$ (i.e. the denominator of $W_k(R(t),p_k(t))$ is canceled with  $p_{k,2}(t)^{\ell}C_{j}(t)^{\ell\deg(q_k)}$). Therefore, from $\num({H}^{\cP{\cal Q}}_{k}(t,R(t)))^\ell=\epsilon^\ell e_k(t)$, and taking into account that $\|\num(W_k(R,p_k))\|\leq \|{H}^{\cP{\cal Q}}_{k}\|^\ell$, we get that
  \[\|\num({H}^{\cP{\cal Q}}_{k}(t,R(t)))\|^\ell=\epsilon^\ell\|e_k\|\leq \epsilon^\ell\|{H}^{\cP{\cal Q}}_{k}\|^\ell,\]
  which implies that  $H^{{\cal P}{\cal Q}}_k(t, R(t))\approx_\epsilon 0$ (see Definition \ref{Def-operator}). Thus,  $S^{{\cal P}{\cal Q}}_\epsilon(t,R(t))\approx_\epsilon 0$, and then ${\cal P}(t)\sim_\epsilon {\cal Q} (R(t))$.\qed
\para

\para

\begin{remark}\label{R-degreex} From the proof of Theorem \ref{Th-resultant}, we have that $\deg_{x_k}(L_k)=\deg_{x_k}(W_k)=\deg(R)$.
\end{remark}

\para

\begin{remark}\label{R-Thresultant} If the tolerance in Theorem \ref{Th-resultant} changes (that is, instead $\epsilon$ we have $\overline{\epsilon}$), Theorem \ref{Th-resultant} holds. More precisely, if
 \[L_k(s, x_k)= (x_kq_{k,2}(s)-q_{k,1}(s))^{\ell}+\overline{\epsilon}^\ell \, W_k(s, x_k),\,\,\,\,\|\num(W_k(R,p_k))\|\leq \|{H}^{\cP{\cal Q}}_{k}\|^\ell,\]
where $\ell=\deg(R)$ and    $\egcd(q_{k,1}, q_{k,2})=1$,
then
 ${\cal
Q}(s)=\left(\frac{q_{1,1}(s)}{q_{1,2}(s)},\,\frac{q_{2,1}(s)}{q_{2,2}(s)}\right)$  is  an  $\overline{\epsilon}$-numerical   reparametrization  of ${\cal P}$.\end{remark}

\para

\para

\para

\begin{center}
\underline{\sf Properties of the $\epsilon$-Numerical Reparametrization ${\cal Q}$}
\end{center}

\para

Let $\cal Q$ be  the  $\epsilon$-numerical   reparametrization  of ${\cal P}$ computed in Theorem \ref{Th-resultant}. In the following, we present some corollaries obtained from Theorem \ref{Th-resultant}, where some properties concerning  $\cal Q$ are provided. In particular, we show that ${\cal Q}$ is   $\epsilon$-proper and we prove that $\deg(\cP)=\deg({\cal Q}) \deg(R).$ This last equality also holds in the symbolic case, and it shows   the expected property that the degree of the rational functions defining the  $\epsilon$-proper parametrization ${\cal Q}$ is lower than the  non  $\epsilon$-proper input parametrization ${\cal P}$.

\para

 \begin{corollary}\label{C-resultant}  ${\cal Q}$ is   $\epsilon$-proper.  \end{corollary}
 \noindent {\sc Proof.} Since  $R(t)=\frac{C_{i}(t)}{C_{j}(t)}\in {\Bbb C}(t)\setminus {\Bbb C}$ is such that   $S^{{\cal P}{\cal P}}_\epsilon(t,s)\approx_\epsilon \num(R(t)-R(s))$, and $\cal Q$ is an  $\epsilon$-numerical   reparametrization  of ${\cal P}$ (see Theorem \ref{Th-resultant}),  from Corollary  \ref{Th-egcd}, we conclude that  ${\cal Q}$ is   $\epsilon$-proper.\qed

\para

\begin{remark}\label{R-Thresultant2}   Corollaries \ref{Th-eindex} and \ref{C-resultant} imply that $\ell=\eindex(\cP)$, where $\ell=\deg(R)$ is introduced in Theorem \ref{Th-resultant}.\end{remark}

\para

 \begin{corollary}\label{C-degree} It holds that
 $\deg(\cP)=\deg({\cal Q}) \deg(R).$  \end{corollary}
\noindent {\sc Proof.} First, we observe that
$\deg_{s}(W_k) = \deg(p_k)$, for $k=1,2$. Indeed, since
\[L_k(s, x_k)=\res_t(G_{k}(t, x_k), s C_{j}(t)-C_{i}(t)),\]
we get that, up to constants in ${\Bbb C}(x_k)\setminus\{0\}$, $$L_k(s, x_k)=\prod_{\{\beta_\ell\,|\,G_k(\beta_{\ell}, x_k)=0\}} sC_{j}(\beta_\ell)-C_{i}(\beta_\ell),$$  (see  Sections 5.8 and 5.9 in \cite{Vander}), and thus $$\deg_{s}(L_k) = \deg_s(s C_{j}(t)-C_{i}(t))\deg_t(G_{k}(t, x_k))= \deg(p_k).$$ Since we are working numerically, we may assume w.l.o.g that $\deg_{s}(W_k)=\deg_{s}(L_k)$. On the other side, from Theorem \ref{Th-resultant}, we have that
\[L_k(s, x_k)= (x_kq_{k,2}(s)-q_{k,1}(s))^{\ell}+\epsilon^{\ell} W_k(s, x_k),\quad k=1,2.\]
Since we are working numerically, we may assume w.l.o.g that $$\deg_{s}(W_k) =\deg_s((x_kq_{k,2}(s)-q_{k,1}(s))^{\ell})= \ell\,\deg(q_k).$$
Therefore,  $\ell\,\deg(q_k)=\deg(p_k),\,k=1,2,$ which implies that
 \[\deg(\cP)=\deg({\cal Q}) \ell=\deg({\cal Q}) \deg(R) \]
 (from Theorem \ref{Th-resultant}, we have that $\ell=\deg(R)$ ).
\qed

\para

 \begin{corollary}\label{C-resultant1} Under the conditions of Theorem \ref{Th-resultant}, it holds that
 \[\res_t(p_{k, 2}(t), s C_{j}(t)-C_{i}(t))= q_{k,2}(s)^{\ell}+\epsilon^\ell b_k(s),\quad b_k\in {\Bbb C}[s],\quad k=1,2.\]
\end{corollary}\noindent {\sc Proof.} From Theorem \ref{Th-resultant} and Corollary \ref{C-degree}, we have that
 \[L_k(s, x_k)= (x_kq_{k,2}(s)-q_{k,1}(s))^{\ell}+\epsilon^\ell W_k(s, x_k),\,\,\,k=1,2,\]
and  $\deg_{s}(W_k)=\ell\,\deg(q_k)=\deg(p_k)$.
Let $L^*_k(s,x_k,x_3)$ be the homogeneous form of the polynomial $L_k(s,x_k)$ w.r.t. the variable  $x_k$. Using the specialization resultant property (see Section 2), we deduce that
\[L^*_k(s, x_k,x_3)=\res_t(x_kp_{k, 2}(t)-x_3 p_{k, 1}(t), s C_{j}(t)-C_{i}(t))=\]\[=(x_kq_{k,2}(s)-x_3q_{k,1}(s))^\ell+\epsilon^\ell b_k(s)x_k^\ell+\epsilon^\ell x_3 U^*_k(s,  x_k,x_3),\]
where  $U^*(s,x_k,x_3)$ denotes the homogeneous form of the polynomial $W_k(s,x_k)-b_k(s)x_k^\ell \in ({\Bbb C}[s])[x_k]$, and $b_k$ is the leading coefficient of $W_k$ w.r.t. $x_k$ that is,  the coefficient of $W_k$ w.r.t. $x_k^{\ell}$ (see  Remark \ref{R-degreex}). Hence, from the specialization resultant property (see Section 2), we get
\[L^*_k(s, 1,0)=\res_t(p_{k, 2}(t), s C_{j}(t)-C_{i}(t))=q_{k,2}(s)^\ell+\epsilon^\ell b_k(s),\,\,\,k=1,2.\]
\qed

 \begin{corollary}\label{C-resultant0} Under the conditions of Theorem \ref{Th-resultant}, it holds that
 \begin{enumerate}\item The rational function $q_k(s):=q_{k,1}(s)/q_{k,2}(s)$ can be obtained by simplifying  the root in the variable $x_k$ of the polynomial $\displaystyle{\frac{\partial^{\ell-1} L_k}{\partial^{\ell-1} x_k}(s, x_k)},\,\,k=1,2$.
 \item The rational function $q_k(s):=q_{k,1}(s)/q_{k,2}(s)$  can be obtained by simplifying  the   rational function $\frac{-\coeff(L_k,x_k^{\ell-1})/\ell}{\coeff(L_k,x_k^\ell)},\,k=1,2$, where $\coeff(pol,var)$ denotes de coefficient of a polynomial $pol$ w.r.t.  $var$.\end{enumerate}
\end{corollary}\noindent {\sc Proof.}  In order to prove statement 1, we  write
 \[L_k(s, x_k)= (x_kq_{k,2}(s)-q_{k,1}(s))^{\ell}+\epsilon^\ell W_k(s, x_k),\,\,\,\,\,k=1,2,\]
where  $\deg_{s}(W_k)=\ell\,\deg(q_k)=\deg(p_k)$ (see Corollary \ref{C-degree}). Thus,
\[\frac{\partial^{\ell-1} L_k}{\partial^{\ell-1} x_k}(s, x_k)=\ell!x_kq_{k,2}(s)^\ell-\ell!q_{k,1}(s)q_{k,2}(s)^{\ell-1}+\ell! \epsilon^\ell x_k b_k(s)-(\ell-1)!\epsilon^\ell a_k(s),\]
where $b_k(s)$ is  the coefficient of $W_k$ w.r.t. $x_k^{\ell}$  (see   Remark \ref{R-degreex}), and  $a_k(s)$ is the  coefficient of $W_k$ w.r.t. $x_k^{\ell-1}$.
The root in the variable $x_k$ of this polynomial is
\[\frac{q_{k,1}(s)q_{k,2}(s)^{\ell-1}+\epsilon^\ell a_k(s)/\ell}{q_{k,2}(s)^\ell+\epsilon^\ell  b_k(s)},\,\,\,\,\,k=1,2.\]
Statement 2 is obtained from statement 1 using that: $$\coeff(L_k,x_k^\ell)=q_{k,2}(s)^\ell+\epsilon^\ell  b_k(s),\,\,\mbox{and}\,\,\, \coeff(L_k,x_k^{\ell-1})=-\ell q_{k,2}(s)^{\ell-1}q_{k,1}(s)-\epsilon^\ell  a_k(s),\,\, k=1,2.\, \mbox{ \qed}$$


\para

\para

\noindent
In the following, we consider the parametrization
\[\widetilde{\cal Q}(s)=\left(\widetilde{q}_{1},\widetilde{q}_{2}\right)=\left(\frac{\widetilde{q}_{1,1}}{\widetilde{q}_{1,2}},\frac{\widetilde{q}_{2,1}}{\widetilde{q}_{1,2}}\right)=\left(\frac{q_{1,1}(s)q_{1,2}(s)^{\ell-1}+\epsilon^\ell a_1(s)/\ell}{q_{2,2}(s)^\ell+\epsilon^\ell  b_1(s)},\,\frac{q_{2,1}(s)q_{2,2}(s)^{\ell-1}+\epsilon^\ell a_2(s)/\ell}{q_{2,2}(s)^\ell+\epsilon^\ell  b_2(s)}\right)\]
%
%
(see Corollary \ref{C-resultant0}). Corollary \ref{C-resultant0} states that  $\widetilde{\cal Q}$ can be further simplified by removing the approximate gcd from the numerator and denominator.  For instance, one may use $\QRGCD$ algorithm to compute an approximate gcd of two univariate polynomials (see more details before Definition \ref{Def-approxindex}). The simplification of $\widetilde{\cal Q}$ provides  the searched rational parametrization ${\cal Q}(s)=\left(\frac{q_{1,1}(s)}{q_{1,2}(s)}, \frac{q_{2,1}(s)}{q_{1,2}(s)}\right)$ (see Corollary \ref{C-resultant0}).






\section{Error Analysis}

In this section, we show the relation between the input curve and the output curve. For this purpose, let   $\cal C$ be the  input curve defined by the parametrization  $\cP=(\frac{p_{1,1}}{p_{1,2}}, \frac{p_{2,1}}{p_{1,2}})$,
with $\index(\cP)=1$. Since ${\cal P}$ is expected to be given with perturbed float coefficients, we may assume w.l.o.g that $\gcd(p_{k,1}, p_{1,2})=1$, $k=1,2$.

\para

 Let  $\widetilde{\cal D}$ be the  curve defined by the parametrization  $\widetilde{\cal Q}=\left(\frac{\widetilde{q}_{1,1}}{\widetilde{q}_{1,2}},\frac{\widetilde{q}_{2,1}}{\widetilde{q}_{1,2}}\right)$ with $\gcd(\widetilde{q}_{k,1}, \widetilde{q}_{1,2})=1$, for   $k=1,2$.    We  may assume w.l.o.g   that $\index(\widetilde{\cal Q})=1$ (note that we are working numerically and then,   with
probability almost one $\deg_t(S)=1$, where $S$ is the polynomial introduced in Section 2).

 \para

\noindent
Observe that from  Corollary \ref{C-resultant1},
 \[\res_t(p_{k, 2}(t), s C_{j}(t)-C_{i}(t))= q_{k,2}(s)^{\ell}+\epsilon^\ell b_k(s)=\widetilde{q}_{k,2}(s),\quad k=1,2\]
which implies that if $p_{1,2}=p_{2,2}$, then  $\widetilde{q}_{1,2}=\widetilde{q}_{2,2}$. In addition, we  may assume w.l.o.g that $\deg(p_{k,1})=\deg(p_{1,2})$ and $\deg(\widetilde{q}_{k,1})=\deg(\widetilde{q}_{1,2})$, for $k=1,2$ (otherwise, one may apply both parametrizations a birational parameter transformation). Thus, $$\deg(p_1)=\deg(p_2)=\deg(p_{1,1})=\deg(p_{2,1})=\deg(p_{1,2}),\,\,\,\,\mbox{ and}\,\,\, $$$$ \deg(\widetilde{q}_1)=\deg(\widetilde{q}_2)=\deg(\widetilde{q}_{1,1})=\deg(\widetilde{q}_{2,1})=\deg(\widetilde{q}_{1,2}).$$

\para

Under these conditions,  in the following theorem we prove that $\deg(f)=\deg(h)$, where $f\in {\Bbb C}[x_1,x_2]$ is the polynomial defining implicitly the curve ${\cal C}$, and $h\in {\Bbb C}[x_1,x_2]$ is the polynomial defining implicitly the curve $\widetilde{\cal D}$.

\para

\begin{theorem} \label{T-degreecurves}  The curves   $\cal C$ and  $\widetilde{\cal D}$ have the same degree. \end{theorem}
\noindent {\sc Proof.} First, taking into account that $$\deg(p_{1,1})=\deg(p_{2,1})=\deg(p_{1,2}),\quad \deg(\widetilde{q}_{1,1})=\deg(\widetilde{q}_{2,1})=\deg(\widetilde{q}_{1,2}),$$ and applying Theorem 6.3.1 in \cite{SWP}, one has that   all the infinity points of both curves are reachable by the corresponding projective parametrizations. In addition, if $(s_0, w_0)$ is such that $\widetilde{q}^*_{1,2}(s_0, w_0)=0$, where  $$\widetilde{\cal Q}^*(s, w)=\left({\widetilde{q}^*_{1,1}(s,w)},\,{\widetilde{q}^*_{2,1}(s,w)}, {\widetilde{q}^*_{1,2}(s,w)}\right)$$ is the projective parametrization of  the
projective curve $\widetilde{\mathcal D}^*$, then  $\widetilde{\cal Q}^*(s_0, w_0)$ generates an infinity point of    $\widetilde{\mathcal D}^*$. Observe that since $\deg(\widetilde{q}_{1,1})=\deg(\widetilde{q}_{2,1})=\deg(\widetilde{q}_{1,2})$, it holds that $\widetilde{q}^*_{1,2}(s_0, w_0)=0$ if and only if $\widetilde{q}_{1,2}(s_0)=0$.  Thus,   the homogeneous form of maximum degree of $h$ is given as $$h_r(x_1, x_2)=\prod_{\{s_i\,|\,\widetilde{q}_{1,2}(s_i)=0\}} (x_2{\widetilde{q}_{1,1}(s_i)}-{\widetilde{q}_{2,1}(s_i)}x_1).$$
Since we are working with approximate mathematical objects, we may assume w.l.o.g that the polynomial $\widetilde{q}_{1,2}$ does not have multiple roots, and then $r=\deg(h)=\deg(\widetilde{q}_{1,2})$. In addition, since
 \[\res_t(p_{1, 2}(t), s C_{j}(t)-C_{i}(t))= \widetilde{q}_{1,2}(s), \]  (see Corollary \ref{C-resultant1}), we get that, up to constants in ${\Bbb C}\setminus\{0\}$, $$\widetilde{q}_{1,2}(s)=\prod_{\{\gamma_\ell\,|\,p_{1, 2}(\gamma_{\ell})=0\}} sC_{j}(\gamma_\ell)-C_{i}(\gamma_\ell),$$  (see  Sections 5.8 and 5.9 in \cite{Vander}), and thus $\deg(\widetilde{q}_{1,2})=\deg_s(s C_{j}(t)-C_{i}(t))\deg(p_{1, 2})=\deg(p_{1, 2})$. Hence,  $r=\deg(h)=\deg({p}_{1,2})$.\\

\noindent
Reasoning similarly with the projective parametrization  of  the
projective input curve ${\mathcal C}^*$, we have that  all the infinity points of ${\mathcal C}^*$ are reachable by the projective parametrization
 $${\cal P}^*(t, w)=\left({p^*_{1,1}(t,w)},\,{p^*_{2,1}(t,w)}, {p^*_{1,2}(t,w)}\right).$$
Hence, the homogeneous form of maximum degree of the polynomial $f$ is given as $$f_d(x_1, x_2)=\prod_{\{t_i\,|\,{p}_{1,2}(t_i)=0\}} (x_2{{p}_{1,1}(t_i)}-{{p}_{2,1}(t_i)}x_1).$$
Since we are working  with approximate mathematical objects, we may assume w.l.o.g that the polynomial ${p}_{1,2}$ does not have multiple roots, and then $d=\deg(f)=\deg({p}_{1,2})$. Therefore,  we conclude that   $d=\deg(f)=\deg(h)=\deg({p}_{k})=\deg(\widetilde{q}_{k})$. \qed

      \para

\begin{remark}\label{R-infinity}   Observe that $x_k,\,k=1,2$ does not divide $f_d$. That is, $f_d(1,0)f_d(0,1)\not=0$. Indeed, if $f_d(1,0)=0$ (similarly if  $f_d(0,1)=0$), one gets that there exists $t_i\in {\Bbb C}$ such that ${p}_{1,2}(t_i)={p}_{2,1}(t_i)=0$. This is impossible, because we have assumed that $\gcd({p}_{1,2},{p}_{2,1})=1$.
Reasoning similarly one proves that $x_k,\,k=1,2$ does not divide $h_d$.
\end{remark}

 \para

From Theorem \ref{T-degreecurves}, we deduce   that the curves ${\mathcal C}$ and $\widetilde{\mathcal D}$ have the same behavior  at
infinity (see Theorem \ref{Th-asymptotes}). More precisely, we show that the homogeneous form of maximum degree of $h$ is equal to the homogeneous form of maximum degree of $f$.

\para

 \begin{theorem}\label{Th-asymptotes}  The implicit equation defining the curves ${\mathcal C}$ and $\widetilde{\mathcal D}$ have the same homogeneous form of
maximum degree, and hence both curves have the same points at
infinity.  \end{theorem}
\noindent {\sc Proof.} First,  by applying Theorem 8 in
\cite{Sen2}, one has that
\[f(x_1,x_2)^{\index(\cP)}=\resultant_t(G_1(t,x_1), G_2(t,x_2)),\, \mbox{where}\quad G_{k}(t, x_k)=x_kp_{1, 2}(t)-p_{k, 1}(t),\,\, k=1,2.\]
We recall that we are assuming that $\index(\cP)=1$.
Now, we consider the polynomials $$L_k(s, x_k)=\res_t(G_{k}(t, x_k), s C_{j}(t)-C_{i}(t)),\,k=1,2,$$ introduced in Theorem \ref{Th-resultant}. By Corollary \ref{C-resultant0}, we have that
$$L_k(s, x_k)=x_k^\ell \widetilde{q}_{1,2}-\ell x_k^{\ell-1} \widetilde{q}_{k,1}+A_k(s,x_k),\quad \deg_{x_k}(A_k)\leq \ell-2,\,\,\,\,k=1,2.$$
Let us prove that there exists a non empty open subset   $\Omega \subset {\Bbb C}^2$, such that for every $q\in {\Omega}$ with $f(q)=0$, it holds that $R(q)=0$, where $R(x_1,x_2):=\resultant_s(L_1, L_2)$. Thus, one would deduce that $f$ divides $R$. Indeed,
first we observe that $R\not=0$ because does not exist any factor depending on $s$ that divides $L_k$ (note that $\gcd(q_{1,2},q_{k,1})=1$). Now,  let
$$\Omega:=\{q\in {\Bbb C}^2\,|\,\lc(G_1,t)(q)\lc(G_2,t)(q) D_2(q) C_{j}(\cP^{-1}(q))\not=0\},$$
where $\cP^{-1}(x_1,x_2)=D_1(x_1,x_2)/D_2(x_1,x_2)$ (note that $\index(\cP)=1$ and then, there exists the inverse of $\cP$ in ${\Bbb C}(x_1,x_2)\setminus{\Bbb C}$).
Observe that $\Omega$ is a non empty open subset of   ${\Bbb C}^2$ since $$\lc(G_1,t)(x_1)\lc(G_2,t)(x_2) D_2(x_1,x_2) C_{j}(\cP^{-1}(x_1,x_2))\not=0.$$ Now, let $q=(x^0_1,x^0_2)\in {\Omega}$ be such that $f(q)=0$ (note that ${\cal C}$ and ${\Bbb C}^2\setminus {\Omega}$ intersect at finitely many points). Since $\lc(G_j,t)(q)\not=0,\,j=1,2$,  by the resultant property (see Section 2), there exists $t_0\in {\Bbb C}$ such that $G_k(t_0,x^0_k)=0,\,k=1,2$. In addition, since $q\in {\Omega}$, one has that there exists  $s_0\in {\Bbb C}$ such that $s_0 C_{j}(t_0)-C_{i}(t_0)=0$ (note that $t_0=\cP^{-1}(q)$, and $C_{j}(t_0)\not=0$). Then,  since $L_k(s, x_k)=\res_t(G_{k}(t, x_k), s C_{j}(t)-C_{i}(t)),$ we get that $L_k(s_0, x^0_k)=0,\,k=1,2.$ Hence, by the specialization of the resultant property (see Section 2), we deduce that $$R(q)=\resultant_s(L_1(s,x_1), L_2(s,x_2))(q)=\resultant_s(L_1(s,x^0_1), L_2(s,x^0_2))=0.$$

\noindent
Thus,
\[R(x_1,x_2)=f(x_1,x_2)m(x_1,x_2),\quad m\in {\Bbb C}[x_1,x_2].\]
Since we are  working  with approximate mathematical objects, we may assume w.l.o.g that  $\deg_{\{x_1,x_2\}}(R)=\deg_s(L_1)\deg_s(L_2)$ (see  Sections 5.8 and 5.9 in \cite{Vander}). Then, if we homogenize the above equation with respect to the variables $x_1$ and $x_2$, we get that
\[R^*(x_1,x_2,x_3):=\resultant_s(L_1^*(s,x_1,x_3), L_2^*(s,x_2,x_3))=F(x_1,x_2,x_3)M(x_1,x_2,x_3),\]
where $F, M  \in {\Bbb C}[x_1,x_2,x_3]$ are the homogenization of $f, m$, respectively, with respect to the variables $x_1$ and $x_2$, and
$$L_k^*(s, x_k,x_3)=x_k^\ell \widetilde{q}_{1,2}-\ell x_k^{\ell-1} x_3 \widetilde{q}_{k,1}+x_3^2A_k(s,x_k,x_3),\quad \deg_{\{x_k,x_3\}}(A_k)=\ell-2,\,\,\,\,k=1,2,$$
 is the homogenization of $L_k$ with respect to  $x_1$ and $x_2$. Observe that $x_3$ does not divide $M$ because $\deg_{\{x_1,x_2\}}(R)=\deg_s(L_1)\deg_s(L_2)$.

\para

\noindent
Now, we consider the system defined by the polynomials $$L_1^*=(x_1^\ell \widetilde{q}_{1,2}+x_3^2A_1(s,x_k,x_3))+x_3(-\ell x_1^{\ell-1}  \widetilde{q}_{1,1}),\,  L_2^*=(x_2^\ell \widetilde{q}_{1,2}+x_3^2A_2(s,x_k,x_3))+x_3(-\ell x_2^{\ell-1} \widetilde{q}_{2,1}).$$
Observe that the two equations are independent. Thus, solving from $L_1^*=0$, we have that $x_3=(x_1^{\ell}q_{1,2}+x_3^2A_1(s,x_1,x_3))/(\ell x_1^{\ell-1} \widetilde{q}_{1,1})$. We   substitute it in $L_2^*$, and we obtain the following equivalent system defined by the polynomials $L_1^*$, and $L^*$, where
\[L^*(s,x_1,x_2,x_3):= \widetilde{q}_{1,2}(s)x_1^{\ell-1}x_2^{\ell-1}(-\widetilde{q}_{2,1}(s)x_1+\widetilde{q}_{1,1}(s)x_2)+x_3^2B(s,x_1,x_2,x_3),\,  B \in {\Bbb C}[s,x_1,x_2,x_3].\]
Thus,
$$R^*(x_1,x_2,x_3)=\resultant_s(L_1^*(s,x_1,x_3), L^*(s,x_2,x_3))=F(x_1,x_2,x_3)M(x_1,x_2,x_3).$$
Using the property of specialization of the resultants, we consider $x_3=0$ in the above equality, and we get that (we remind that $x_3$ does not divide $M$)
 \[\resultant_s(\widetilde{q}_{1,2}(s), x_1^{\ell-1}x_2^{\ell-1}(-\widetilde{q}_{2,1}(s)x_1+\widetilde{q}_{1,1}(s)x_2))=\]\[x_1^{\deg(\widetilde{q}_{1,2})(\ell-1)}x_2^{\deg(\widetilde{q}_{1,2})(\ell-1)}\resultant_s(\widetilde{q}_{1,2}(s), (-\widetilde{q}_{2,1}(s)x_1+\widetilde{q}_{1,1}(s)x_2))=f_d(x_1,x_2)m_{\ell}(x_1,x_2),\] where $f_d, m_\ell$ are the  homogeneous form of
maximum degree of $F, M$, respectively.

\para

\noindent
On the other side,   by applying Theorem 8 in
\cite{Sen2}, one also has that
\[h(x_1,x_2)^{\index(\widetilde{\cal Q})}=\resultant_t( \widetilde{G}_1(t,x_1),  \widetilde{G}_2(t,x_2)),\quad\mbox{where}\quad \widetilde{G}_{k}(t, x_k)=x_k\widetilde{q}_{1, 2}(t)-\widetilde{q}_{k, 1}(t),\,\,k=1,2.\]
 We recall that $\index(\widetilde{\cal Q})=1$.
  Since we are  working  with approximate mathematical objects, similarly as above we may assume  that $\deg_{\{x_1,x_2\}}(h)=\deg_t(\widetilde{G}_1)\deg_t(\widetilde{G}_2)$. Then, if we homogenize the above equation with respect to the variables $x_1$ and $x_2$, we get that
  \[H(x_1,x_2,x_3)=\resultant_t( \widetilde{G}^*_1(t,x_1,x_3),  \widetilde{G}^*_2(t,x_2,x_3)),\quad \mbox{where}\,\,\, \widetilde{G}^*_{k}(t, x_k,x_3)=x_k\widetilde{q}_{1, 2}(t)-\widetilde{q}_{k, 1}(t)x_3,\]
and $H$ is the homogenization of $h$ with respect to the variables $x_1$ and $x_2$.  Observe that $x_3$ does not divide $H$ because $\deg_{\{x_1,x_2\}}(h)=\deg_t(\widetilde{G}_1)\deg_t(\widetilde{G}_1)$.\\
Now, reasoning as above, we have that the  system defined by the polynomials $ \widetilde{G}^*_{1}$ and $\widetilde{G}^*_{2}$ is equivalent to the system defined by $\widetilde{G}^*_1$ and the polynomial $\widetilde{G}^*=-\widetilde{q}_{2,1}(s)x_1+\widetilde{q}_{1,1}(s)x_2$.
Thus,   \[H(x_1,x_2,x_3)=\resultant_t(\widetilde{G}^*_1(t,x_1,x_3),  \widetilde{G}^*(t,x_1,x_2)).\]
  Using the property of specialization of the resultants, we consider $x_3=0$ in the above equality, and we get that (observe that $x_3$ does not divide $H$)
\[\resultant_s(\widetilde{q}_{1,2}(s), -\widetilde{q}_{2,1}(s)x_1+\widetilde{q}_{1,1}(s)x_2)=h_d(x_1,x_2),\]
where $h_d$ is the  homogeneous form of
maximum degree of $H$ (we recall that $d=\deg(f)=\deg(h)$, see Theorem \ref{T-degreecurves}).  Thus, since
\[f_d(x_1,x_2)m(x_1,x_2)=x_1^{\deg(\widetilde{q}_{1,2})(\ell-1)}x_2^{\deg(\widetilde{q}_{1,2})(\ell-1)}\resultant_s(\widetilde{q}_{1,2}(s), (-\widetilde{q}_{2,1}(s)x_1+\widetilde{q}_{1,1}(s)x_2))=\]\[x_1^{\deg(\widetilde{q}_{1,2})(\ell-1)}x_2^{\deg(\widetilde{q}_{1,2})(\ell-1)}h_d(x_1,x_2),\]
  and  $x_k,\,\,k=1,2,$ does not divide $f_d$ (see Remark \ref{R-infinity}),
we conclude that $h_d=f_d$. Hence, ${\mathcal C}$ and $\widetilde{\mathcal D}$ have the same homogeneous form of
maximum degree, and then both curves have the same degree and the same points at
infinity.  \qed

\para

\para

\para

 As we stated above,  the parametrization  $\widetilde{\cal Q}$ should be further simplified to obtain the searched   parametrization $\cal Q$ (note that by Theorem \ref{T-degreecurves},   $\deg(\widetilde{q}_{k})=\deg(p_k)$, and from Corollary \ref{C-degree}, we look for $\cal Q$ such that $\ell\,\deg(q_k)=\deg(p_k)$). However, when we simplify  $\widetilde{\cal Q}$, the curve $\widetilde{\cal D}$ defined by $\widetilde{\cal Q}$  changes (the infinity points are not the same because the numerical simplification). This is the expected situation because that in fact,  the input parametrization $\cP$ and the output parametrization $\cal Q$ have different degrees (see Corollary \ref{C-degree}).

\para

 In order to analyze the behavior at affine points, we study the closeness of the curves $\cal C$ and ${\cal D}$,  where $\cal D$ is the curve defined by the simplified parametrization   ${\cal Q}=\left(\frac{{q}_{1,1}}{{q}_{1,2}},\frac{{q}_{2,1}}{{q}_{2,2}}\right)$ (note that by Corollary 3,  $\eindex(Q)=1$),  and $\cal C$ is the curve defined by  $\cP=(\frac{p_{1,1}}{p_{1,2}}, \frac{p_{2,1}}{p_{2,2}})$.  For this purpose, we first assume that  $\deg(p_{i,1})=\deg(p_{i,2})$, and $\deg({q}_{i,1})=\deg({q}_{i,2})$,\, $i=1,2$ (otherwise, one applies both parametrizations a birational parameter transformation).  In addition, let $\|p\|:=\max\{\|p_{1,1}\|, \|p_{2,1}\|, \|p_{1,2}\|, \|p_{2,2}\|\}$, and $\|q\|:=\max\{\|q_{1,1}\|, \|q_{2,1}\|, \|q_{1,2}\|, \|q_{2,2}\|\}$.

\para

Finally, we also assume that Theorem \ref{Th-resultant} holds and then, $\cal Q$ is an $\epsilon$-proper reparametrization of $\cal P$ (see Corollary \ref{C-resultant}). If Theorem \ref{Th-resultant} does not hold, one applies Remark \ref{R-Thresultant}, and then $\cal Q$ an $\overline{\epsilon}$-proper reparametrization of $\cal P$. In this case, the formula obtained in Theorem \ref{Th-error} is the same but it involves  $\overline{\epsilon}$ instead $\epsilon$.

\para

Under these conditions,  in order to analyze the behavior at affine points, we restrict us to an interval where the parametrizations $\cP $ and ${\cal Q}$ are well defined. Thus,   the general strategy we  follow is to show that almost any affine
real point on $\cal D$ is at small distance of an affine real point on $\cal C$ (and reciprocally).

\para

For this purpose,  we consider the interval   $I:=(d_1, d_2) \subset {\Bbb R}$ satisfying that for all  $t_0 \in I$, there exists $M\in {\Bbb N}$ such that $|{q}_{i, 2}(R(t_0))|\geq M$, and $|p_{i, 2}(t_0)|\geq M$,\,$i=1,2$. Note that we can decompose $\Bbb R$ as union of finitely many intervals, $I_j, j=1,\ldots,n$, satisfying the above condition (that is, we consider intervals when no roots of the denominators of the parametrizations appear; see \cite{PSV}). Then, we may reason similarly as we do in Theorem \ref{Th-error} for
each interval $I_j, j=1,\ldots,n$.

\para

\begin{theorem} \label{Th-error} The following statements hold:
\begin{enumerate}
\item  Let  $I:=(d_1, d_2) \subset {\Bbb R}$, and $M\in {\Bbb N}$ such that for every $t_0 \in I$, it holds that
$ |{q}_{i, 2}(R(t_0))|\geq M$, and $|p_{i, 2}(t_0)|\geq M$  for $i=1,2$.  Let   $d:=\max\{|d_1|, |d_2|\}$.   Then,  for every  $t_0 \in I$,
 $$|p_i(t_0)-{q}_i(R(t_0))| \leq 2/M^2 \epsilon   \,C  \|p\|\|q\|,\quad i=1,2,$$ where $$C=\left\{\begin{array}{l} \displaystyle{\frac{d^{\deg(\cP)+1}}{(d-1)^{1/\ell}}}\,\qquad\quad \mbox{if} \,\, d>1,\\
    \displaystyle{\frac{1}{(1-d)^{1/\ell}}}\,\qquad\quad  \mbox{if} \,\, d<1,\\
    \\
    \ell^{1/\ell} \deg(\cP)^{1/\ell} \qquad  \mbox{if}\,\, d=1.\end{array}\right. $$
     \item  ${\mathcal C}_{t\in I}$ is contained in the offset region of ${\mathcal D}_{s\in J}$ at distance $4\sqrt{2}/M^2 \epsilon   \,C\, \|p\|\|q\|$, where $J=R(I).$
        \item  ${\mathcal D}_{s\in J}$ is contained in the offset region of ${\mathcal C}_{t\in I}$ at distance $4\sqrt{2}/M^2 \epsilon   \,C\, \|p\|\|q\|$, where $J=R(I).$
       \end{enumerate}
\end{theorem}
\noindent {\sc Proof.} First, we observe that statement (1) implies statements (2) and (3). For this purpose, we note that for  almost all affine real point $Q \in {\mathcal D}$ there exists an affine real point ${P} \in {{\mathcal C}}$ such that
$$ \|P-Q\|_2 \leq 2\sqrt{2}/M^2 \epsilon   \,C\, \|p\|\|q\|.$$ Indeed, using statement (1), we have that
$$ \|P-Q\|_2=\sqrt{ (p_1(t_0)-{q}_1(R(t_0)))^2+(p_2(t_0)-{q}_2(R(t_0))^2}\leq$$$$
\sqrt{ (2/M^2 \epsilon   \,C \|p\|\|q\|)^2 +(2/M^2 \epsilon   \,C  \|p\|\|q\|)^2} \leq    2\sqrt{2}/M^2 \epsilon   \,C\, \|p\|\|q\|.$$
 Now, reasoning as in Section 2.2 in \cite{Farouki}, we deduce    statements (2) and (3).\\

\noindent
Now, we prove statement (1). For this purpose, from the proof of Theorem \ref{Th-resultant}, we have that  $${H}^{\cP{\cal Q}}_{i}(t,R(t))^\ell=(p_{i, 1}(t){q}_{1, 2}(R(t))-{q}_{i, 1}(R(t))p_{1, 2}(t))^\ell=\epsilon^\ell e_i(t),\quad \mbox{where}$$
\[e_i(t)=-\num(W_i(R(t),p_i(t)))=e_{i,0}+e_{i,1}t+\ldots+ e_{i,n_{i}}t^{n_i}\in{\Bbb C}[t],\quad \mbox{and}\quad\]\[ \|e_i\|=\|\num(W_i(R,p_i))\|\leq \|{H}^{\cP{\cal Q}}_{i}\|^\ell.\] In addition, since $e_i(t)=-\num(W_i(R(t),p_i(t)))$, we have that $n_i:=\deg(e_i)\leq \ell \deg(\cP)$ for $i=1,2$. Indeed, since $\deg_{x_i}(W_i)=\ell$ (see Remark \ref{R-degreex}), we deduce that
 $$\deg(e_i)\leq \mbox{max}\{\deg(R)\deg_{t}(W_i), \ell\deg({\cal P})\}\leq \mbox{max}\{\ell \deg({\cal P}), \ell\deg({\cal P})\} = \ell \deg(\cP).$$
In these conditions, for every $t_0 \in I$, if $d\not=1$, it holds that $$|{H}^{\cP{\cal Q}}_{i}(t_0,R(t_0))^\ell|=\epsilon^\ell |e_i(t_0)|\leq \epsilon^\ell\|{H}^{\cP{\cal Q}}_{i}\|^\ell(|e_{i,0}|+|e_{i,1}||t_0|+\ldots+ |e_{i,n_{i}}||t_0|^{n_i})\leq$$

\begin{equation}\label{eq4}
  \epsilon^\ell \|{H}^{\cP{\cal Q}}_{i}\|^\ell (1+d+\ldots +d^{n_i})=\epsilon^\ell \|{H}^{\cP{\cal Q}}_{i}\|^\ell\, \frac{d^{n_i+1}-1}{d-1},\quad i=1,2.
\end{equation}
 If $d=1$, then
\begin{equation}\label{eq5}|{H}^{\cP{\cal Q}}_{i}(t_0,R(t_0))^\ell|\leq \epsilon^\ell \|{H}^{\cP{\cal Q}}_{i}\|^\ell (1+|t_0|+\ldots+ |t_0|^{n_i})\leq\epsilon^\ell \|{H}^{\cP{\cal Q}}_{i}\|^\ell (1+1+\ldots+ 1)=\epsilon^\ell \|{H}^{\cP{\cal Q}}_{i}\|^\ell\, {n_i}.\end{equation}
Therefore, we conclude that:
\begin{enumerate}\item If $d>1$, from~\eqref{eq4}, and taking into account that $ |{q}_{i, 2}(R(t_0))|\geq M$, and $|p_{i, 2}(t_0)|\geq M$  for $i=1,2$, we obtain that $$|p_i(t_0)-{q}_i(R(t_0))| =\frac{|{H}^{\cP{\cal Q}}_{i}(t_0,R(t_0))|}{|{q}_{i, 2}(R(t_0))p_{i, 2}(t_0)|}\leq 1/M^2 \epsilon  \|{H}^{\cP{\cal Q}}_{i}\|  \frac{ d^{\deg(\cP)+1/\ell}}{(d-1)^{1/\ell}}\leq $$$$1/M^2 \,\epsilon\, \|{H}^{\cP{\cal Q}}_{i}\|\, \frac{ d^{\deg(\cP)+1}}{(d-1)^{1/\ell}}.$$
 \item If $d<1$, from~\eqref{eq4}, and taking into account that $1-d^{n_i+1}<1$, and    $ |{q}_{i, 2}(R(t_0))|\geq M$, and $|p_{i, 2}(t_0)|\geq M$  for $i=1,2$, we obtain that    $$|p_i(t_0)-{q}_i(R(t_0))| =\frac{|{H}^{\cP{\cal Q}}_{i}(t_0,R(t_0))|}{|{q}_{i, 2}(R(t_0))p_{i, 2}(t_0)|}\leq1/M^2 \,\epsilon \, \|{H}^{\cP{\cal Q}}_{i}\| \frac{ 1}{(1-d)^{1/\ell}}.$$
 \item If $d=1$,  from~\eqref{eq5}, and taking into account that $ |{q}_{i, 2}(R(t_0))|\geq M$, and $|p_{i, 2}(t_0)|\geq M$  for $i=1,2$, we obtain that
 $$|p_i(t_0)-{q}_i(R(t_0))| =\frac{|{H}^{\cP{\cal Q}}_{i}(t_0,R(t_0))|}{|{q}_{i, 2}(R(t_0))p_{i, 2}(t_0)|}\leq1/M^2 \,\epsilon\,  \|{H}^{\cP{\cal Q}}_{i}\|(\ell {\deg(\cP)})^{1/\ell}. $$
\end{enumerate}
Finally, we note that
\[\|{H}^{\cP{\cal Q}}_{i}\|=\|p_{i, 1}(t)q_{i, 2}(s)-q_{i, 1}(s)p_{i, 2}(t)\|\leq 2\|p\|\|q\|.\]\qed

\para

\noindent
From Theorem \ref{Th-error}, we deduce the following corollary: \para

\begin{corollary} Under the conditions of  Theorem \ref{Th-error}, it holds that:
\begin{enumerate}\item If $d\geq 2$, then $C\leq  d^{\deg(\cP)+1}$. \item If  $1< d<  2$, then
$C \leq 2^{\deg(\cP)+1}.$
\end{enumerate}\end{corollary} \noindent {\sc Proof.}
If $d\geq 2$, we use Theorem \ref{Th-error}, and we get that
 $C\leq \frac{ d^{\deg(\cP)+1}}{(d-1)^{1/\ell}}\leq d^{\deg(\cP)+1}$. \\ If  $1< d<  2$,
$$|{H}^{\cP{\cal Q}}_{i}(t_0,R(t_0))^\ell|\leq \epsilon^\ell  \|{H}^{\cP{\cal Q}}_{i}\|^\ell (1+|t_0|+\ldots+ |t_0|^{n_i})\leq $$$$\epsilon^\ell  \|{H}^{\cP{\cal Q}}_{i}\|^\ell (1+2+\ldots+ 2^{n_i})\leq  \epsilon^\ell\,  \|{H}^{\cP{\cal Q}}_{i}\|^\ell 2^{{n_i}+1},\quad i=1,2.$$
Thus,
$$|p_i(t_0)-q_i(R(t_0))| =\frac{|{H}^{\cP{\cal Q}}_{i}(t_0,R(t_0))|}{|q_{1, 2}(R(t_0))p_{1, 2}(t_0)|}\leq 1/M^2 \,\epsilon  \|{H}^{\cP{\cal Q}}_{i}\|  \, 2^{\deg(\cP)+1}\leq  2/M^2 \,\epsilon  \|p\|\|q\|  \, 2^{\deg(\cP)+1}. \,\qed$$


\section{Numeric Algorithm of  Reparametrization for Curves}

In this section, we apply the  results obtained in Section 3 to
derive an algorithm that computes an $\epsilon$-proper
reparametrization of a given approximate improper parametrization of a  plane
curve. We outline this approach, and we illustrate it with some
examples where we also show the error bound obtained applying results in Section 4.

\noindent
\begin{center}
\fbox{\hspace*{2 mm}\parbox{15.7 cm}{ \vspace*{2 mm} {\bf Numeric Algorithm
{\sf  Reparametrization for Curves}.}

 \vspace*{1 mm}

\noindent
{\sc Input}: a tolerance $\epsilon>0$, and a rational parametrization ${\mathcal
P}(t)=\left(\frac{p_{1, 1}(t)}{p_{1, 2}(t)},\,\frac{p_{2, 1}(t)}{p_{2,
2}(t)}\right)\in {\Bbb C}(t)^2,$   $\egcd(p_{i, 1},p_{i, 2})=1,$\,$i=1,2$,  of an algebraic plane curve $\cal C$.\\
{\sc Output}: a  rational parametrization ${\cal Q}(t)=\left(\frac{q_{1, 1}(t)}{q_{1, 2}(t)},\,\frac{q_{2, 1}(t)}{q_{2,
2}(t)}\right)\in {\Bbb C}(t)^2,$ $\egcd(q_{i, 1},q_{i, 2})=1,\,i=1,2$,  such that $\eindex({\cal Q})=1$  and   ${\cal
P}\sim_\epsilon  {\cal Q}(R)$, where  ${R}(t)\in {\Bbb C}(t)\setminus{\Bbb C}$.
\begin{enumerate}
\item[1.] Compute the polynomials  $H_{k}^{\cP\cP}(t,s)=p_{k, 1}(t)p_{k, 2}(s)-p_{k, 1}(s)p_{k, 2}(t),\,\,
k=1,2.$

\item[2.] Compute  $$S^{{\cal P}{\cal P}}_\epsilon(t, s)=\egcd(H_1^{\cP\cP}(t,s), H_2^{\cP\cP}(t,s))\approx_\epsilon
C_m(t)s^m+\cdots+C_0(t),$$  and   $\eindex({\cal P}):=\deg_t(S^{{\cal P}{\cal P}}_\epsilon)$ (see Definition \ref{Def-approxindex}).

\item[3.]  If $\eindex({\cal P})=1$, {\sc Return}   ${\cal Q}(t)={\mathcal P}(t)$, and
$R(t)=t$. Otherwise go to Step 4.

\item[4.]  Consider
$R(t) = \frac{C_{i}(t)}{C_{j}(t)} \in {\Bbb C}(t),$ such that
$C_{j}(t),\, C_{i}(t)$   are two of the polynomials obtained in Step 2 satisfying that $C_jC_i\not\in {\Bbb C}$, $\egcd(C_j, C_i)=1$,   and  $S^{{\cal P}{\cal P}}_\epsilon(t,s) \approx_\epsilon \num(R(t)-R(s)).$

\item[5.]  For $k=1,2$, compute the polynomials (see Theorem \ref{Th-resultant})
\[L_k(s, x_k)=\res_t(G_{k}(t, x_k), s C_{j}(t)-C_{i}(t)), \,\,\mbox{where}\,\,G_{k}(t, x_k)=x_kp_{k, 2}(t)-p_{k, 1}(t).\]

\item[6.]
For $k=1,2$, compute the root in the variable $x_k$ of the polynomial $\displaystyle{\frac{\partial^{\ell-1} L_k}{\partial^{\ell-1} x_k}(s, x_k)}$ (see  Corollary \ref{C-resultant0}), where $\ell:=\deg(R)=\eindex({\cal P})$ (see   Remark \ref{R-Thresultant2}). Let $\widetilde{q}_k(t)=\widetilde{q}_{k,1}(t)/\widetilde{q}_{k,2}(t)$ be this root, and let
$\widetilde{\cal
Q}(t)=\left({\widetilde{q}_{1,1}(t)}/{\widetilde{q}_{1,2}(t)},\,{\widetilde{q}_{2,1}(t)}/{\widetilde{q}_{2,2}(t)}\right)\in {\Bbb C}(t)^2$.

\item[7.] Simplify $\widetilde{\cal
Q}(t)$ (see Remark \ref{R-simplification1}). Let $${\cal
Q}(t)=\left(\frac{q_{1,1}(t)}{q_{1,2}(t)},\,\frac{q_{2,1}(t)}{q_{2,2}(t)}\right) \in {\Bbb C}(t)^2,\,\quad \egcd(q_{k,1},q_{k,2})=1,\,\,k=1,2,$$ be the simplified parametrization. Check whether the following equality holds \[L_k(s, x_k)= (x_kq_{k,2}(s)-q_{k,1}(s))^{\ell}+\epsilon^{\ell} W_k(s, x_k),\,\,\,\,\,\,\,\,\,\|\num(W_k(R,p_k))\|\leq \|{H}^{\cP{\cal Q}}_{k}\|^\ell\] (see Theorem \ref{Th-resultant}). If it does not hold, use Remark \ref{R-Thresultant} and compute $\overline{\epsilon}$.

\item[8.] {\sc Return}   ${\cal Q},$ $R$, and   the message ``{\tt ${\cal
Q}$ is an $\epsilon$-proper reparametrization of $\cal P$}" (or ``{\tt ${\cal
Q}$ is an $\overline{\epsilon}$-proper reparametrization of $\cal P$}", if  Remark \ref{R-Thresultant} is applied).
\end{enumerate}}\hspace{2 mm}}
\end{center}

\begin{remark}\label{R-simplification1} In the $\egcd$ computation, one may use the {\it SNAP} package included in {\tt Maple}. This package is based on \cite{Beckermann1}, \cite{Beckermann2}, \cite{Corless} and \cite{Karmarkar}.
For simplification of $\widetilde{\cal Q}$ in Step 7, we remove an $\egcd$ from its numerator and denominator  under the given tolerance $\epsilon$ (see more details before Definition \ref{Def-approxindex}).
 \end{remark}

\para

In the
following, we illustrate  {\sf Numeric Algorithm Reparametrization
for Curves} with two examples in detail  where  we explain how the algorithm is performed, and we summarize some other  examples in different tables. In these tables   we show the parametrization $\cal P$ defining the input curve $\cal C$, the tolerance $\epsilon$ considered, the
output parametrization $\cal Q$ defining the output curve $\cal D$, the error bound obtained applying Theorem 5, and a figure representing $\cal C$ and $\cal D$. The computations are done with  the computer algebra system {\sf Maple}, and the  number of digits  we are using when calculating with software floating-point numbers is $10$.

\para

\begin{example}   Let $\epsilon=0.01$, and the rational curve $\cal C$  defined by the parametrization
\[{\cal P}(t)=\left(\frac{p_{1,1}(t)}{p_{1,2}(t)},\,
\frac{p_{2,1}(t)}{p_{2,2}(t)}\right)=\]
\[\left(
{\displaystyle \frac {t^{4} - .2502500000 + .0005000000000\,t}{t
^{4} + .2500000000 + .0002500000000\,t^{2}}} , \,{\displaystyle
\frac {t^{2} - .0002500000000}{t^{4} + .2500000000 +
.0002500000000\,t^{2}}}\right)
\]
It is quartic curve and approximately multiple conic curve~(see Figure~\ref{Ex2fig1}).
Using the
{\it SNAP}   package included in Maple, one has that  $\egcd(p_{j, 1}, p_{j, 2})=1,\,\,j=1,2$. We apply the algorithm  and  in Step 1, we compute the polynomials\\

\noindent $H_{1}^{\cP\cP}(t, s)=8004000\,s^{4} + 4000\,s^{4}\,t^{2} - 8004000\,t^{
4} - 1001\,t^{2} + 8000\,s\,t^{4} + 2000\,s + 2\,s\,t^{2} \\
\mbox{} - 4000\,t^{4}\,s^{2} + 1001\,s^{2} - 8000\,t\,s^{4} -
2000\,t - 2\,t\,s^{2},$\\

\noindent $H_{2}^{\cP\cP}(t,s)= 16000000\,t^{4}\,s^{2} + 4000001\,s^{2} - 4000\,t
^{4} - 4000001\,t^{2} - 16000000\,s^{4}\,t^{2} + 4000\,s^{4}.$\\

\noindent Now, we compute the polynomial   $S^{{\cal P}{\cal P}}_\epsilon$. We have,
$$S^{{\cal P}{\cal P}}_\epsilon(t, s)\approx_\epsilon C_0(t)+C_1(t)s+C_2(t)s^2,$$
where $C_0(t)=52160t^2+83t$,\,\,\,$ C_1(t)=-83t-83$,\, and $C_2(t)=-52077.$ Then,  $\eindex({\cal P})=\deg_t(S^{{\cal P}{\cal P}}_\epsilon)=2$ (see Definition \ref{Def-approxindex}).  Now, we   consider the rational function
\[R(t)= \frac{C_0(t)}{C_2(t)}={\frac {52160t^2+83t }{-52077}},\]
(see Step 4), and we determine the polynomials\\

\noindent
$L_1(s,x_1)=\res_t(G_{1}(t, x_1), sC_1(t)-C_0(t))= .06216533631\,x_1 - .2494530177\,s^{2} +\\
.1389703128\,10^{-5}\,s - .6226547971\,10^{-4}s\,x_1 + .2480131435\,10^{-3}s^{2}\,x_1 + .03105156714\,x_1^{2} +\\
.2492050357\,s^{2}\,x_1^{2}  - .6346614958\,10^{-4}s\,x_1^{2} + .5000000000\,s^{4}\,x_1^{2} - 1.\,s^{4}\,x_1 + .2496021856\,10^{-3}s^{3}\,x_1 -.2496021856\,10^{-3}s^{3}\,x_1^{2} + .5000000000\,s^{4} +
.03111380032,$\\

\noindent
$L_2(s,x_2)=\res_t(G_{2}(t, x_2), sC_1(t)-C_0(t))=.6188610482\,10^{-4}x_2 + .2492049406\,x_2\,s +\\
.4992043712\,s^{2} + .2492050042\,10^{-3}s - .1268068308\,10^{-5}\,x_2\,s^{2} + .03110105693\,x_2^{2} -\\
.6356730153\,10^{-4}x_2^{2}\,s  + .2496022171\,x_2^{2}\,s^{2} - .2500000000\,10^{-3}x_2^{2}\,s^{3
} + x_2\,s^{3} + .5007968969\,x_2^{2}\,s^{4}+ .3078504788\,10^{-7},
$\\

\noindent
where $G_{k}(t, x_k)=x_kp_{k, 2}(t)-p_{k, 1}(t),\,\,k=1,2$ (see Step 5).  In Step 6,  we compute the root in the variable $x_k$ of the polynomial $\frac{\partial L_k}{\partial x_k}(s, x_k)$ (see Corollary \ref{C-resultant0}). We get the curve $\widetilde{\mathcal D}$  defined by the rational parametrization (see Figure~\ref{Ex2fig1})
  $$\widetilde{\cal Q}(t)= \left(\frac{-.06216533631+.00006226547971t-.0002480131435t^2+t^4-.0002496021856t^3}{-.0001269322992t+t^4+.06210313427+.4984100713t^2-.0004992043712t^3
},\,\right.$$$$ \left.\frac{-.00006178762808-.2488083914t+.000001266050484t^2-.9984087423t^3
}{-.0001269322992t+t^4+.06210313427+.4984100713t^2-.0004992043712t^3
}\right).$$

\begin{figure}[h]
\begin{center}
\hbox{
\centerline{\epsfig{figure=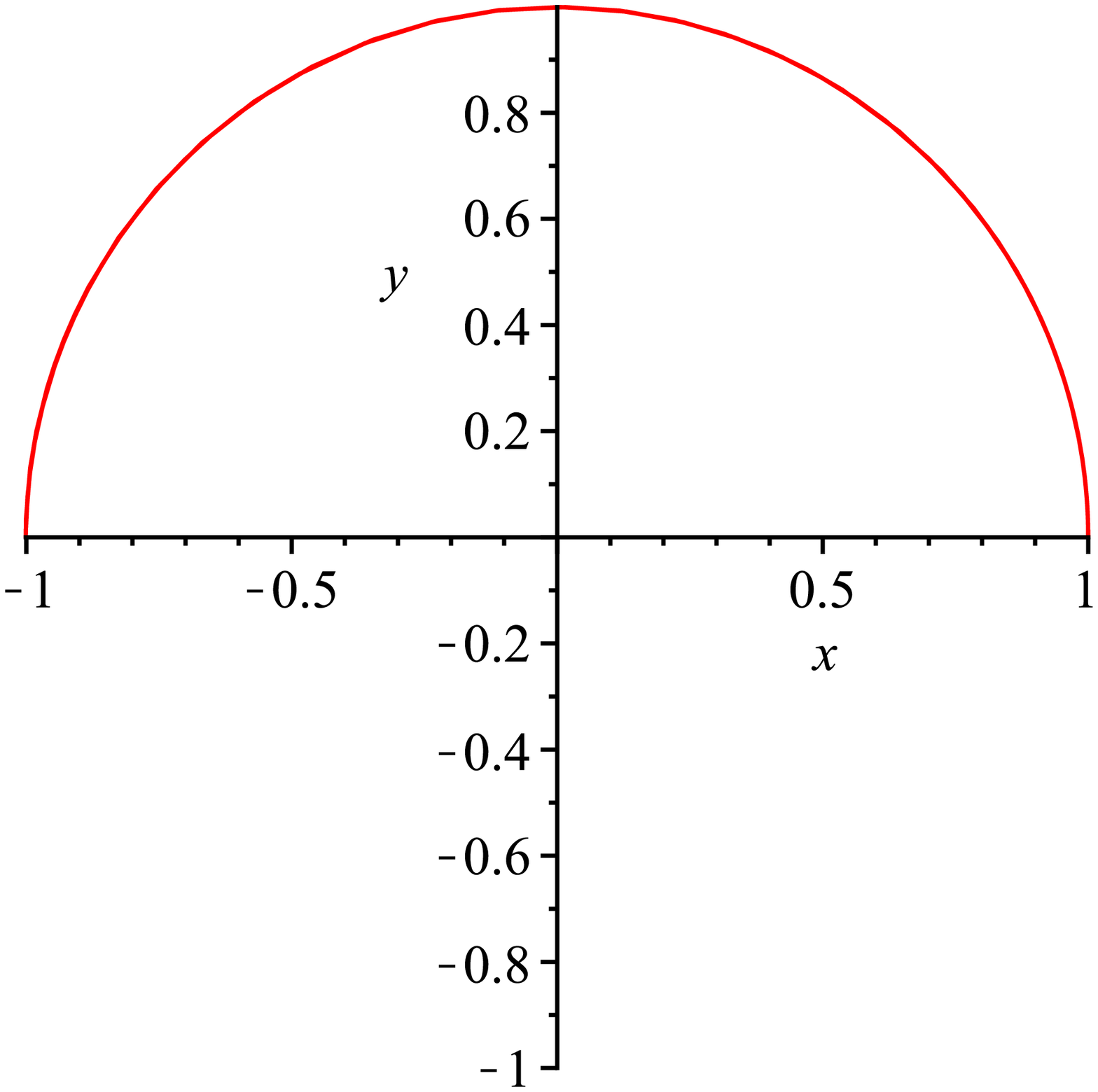,width=5cm} \epsfig{figure=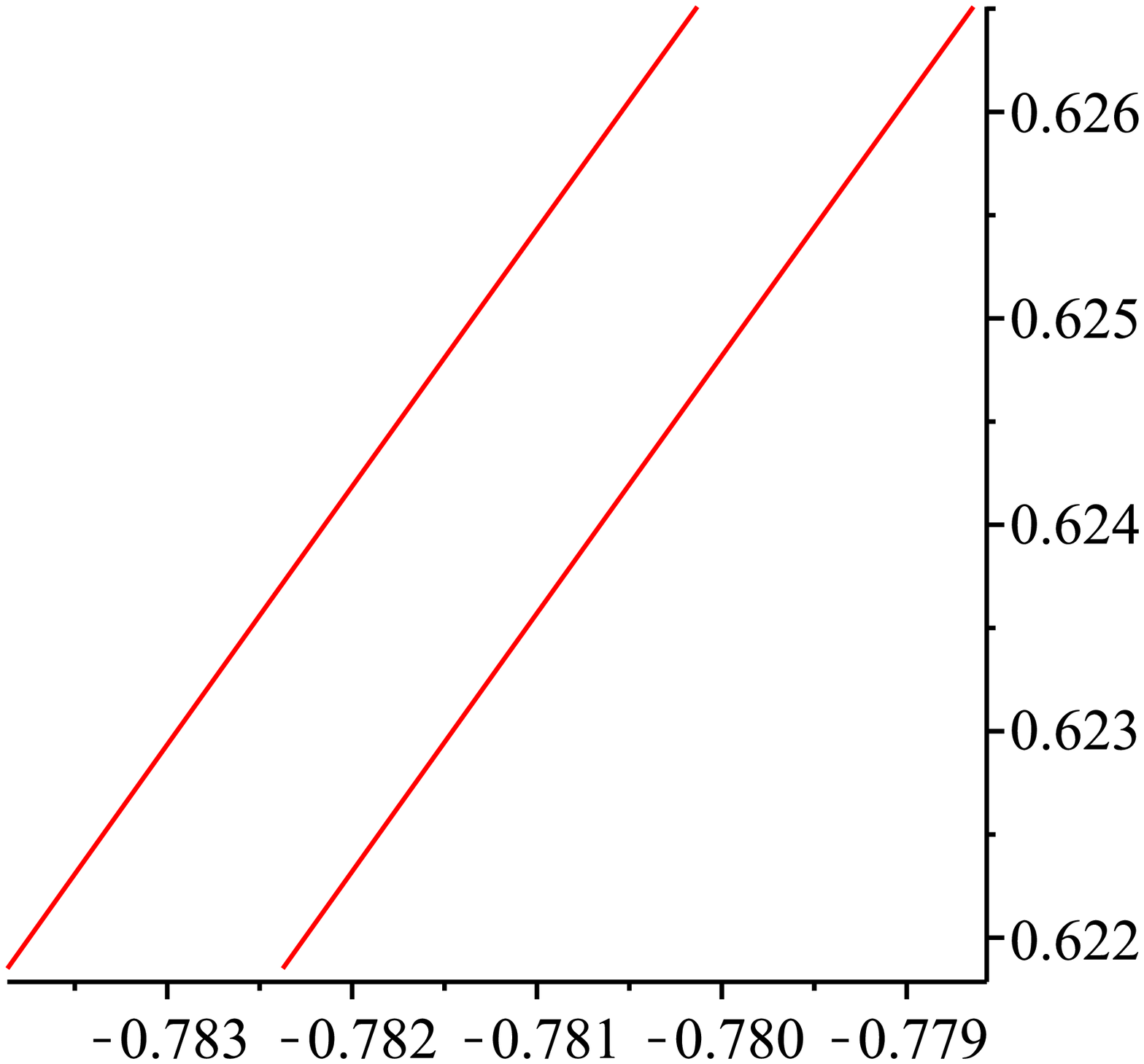,width=5cm,height=5cm} \epsfig{figure=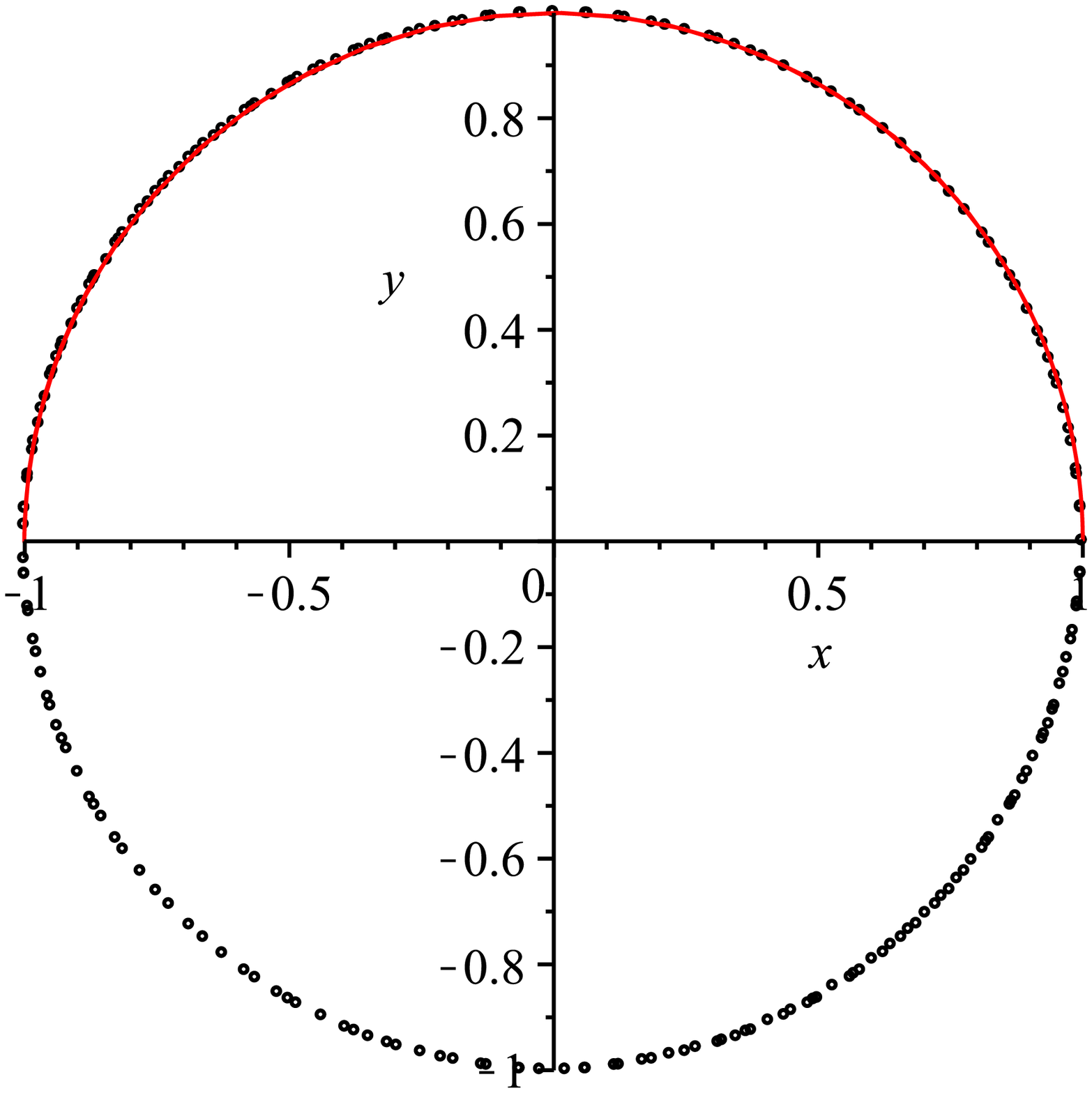,width=5cm,height=5cm}
}}
\end{center}
\vspace*{-1.5cm}
\caption{Input curve $\mathcal C$ (left),  partial enlarged view of $\mathcal C$ (center),  curves $\mathcal C$ and $\widetilde{\mathcal D}$ (right)}
 \label{Ex2fig1}
\end{figure}

\para

Finally,  we simplify $\widetilde{\cal Q}$ by removing certain $\epsilon$-gcds (see Remark \ref{R-simplification1}). We get the curve $\cal D$ defined by the $\epsilon$-numerical reparametrization (see Figure~\ref{Ex2fig2}):
  $$ {\cal Q}(t)= \left(\frac{{q}_{1,1}(t)}{{q}_{1,2}(t)},\,\frac{{q}_{2,1}(t)}{{q}_{2,2}(t)}\right)= $$$$ \left(\frac {t^2+.000005006649227t-.2494538109}{t^2-.0002445955365t+.2492042101}
,\,\frac{-.9984087427t-.0002529376363}{t^2-.0002445955365t+.2492042101}
\right).$$

\begin{figure}[h]
\begin{center}
\hbox{
\centerline{\epsfig{figure=ex2b1.eps,width=5cm,height=5cm} \epsfig{figure=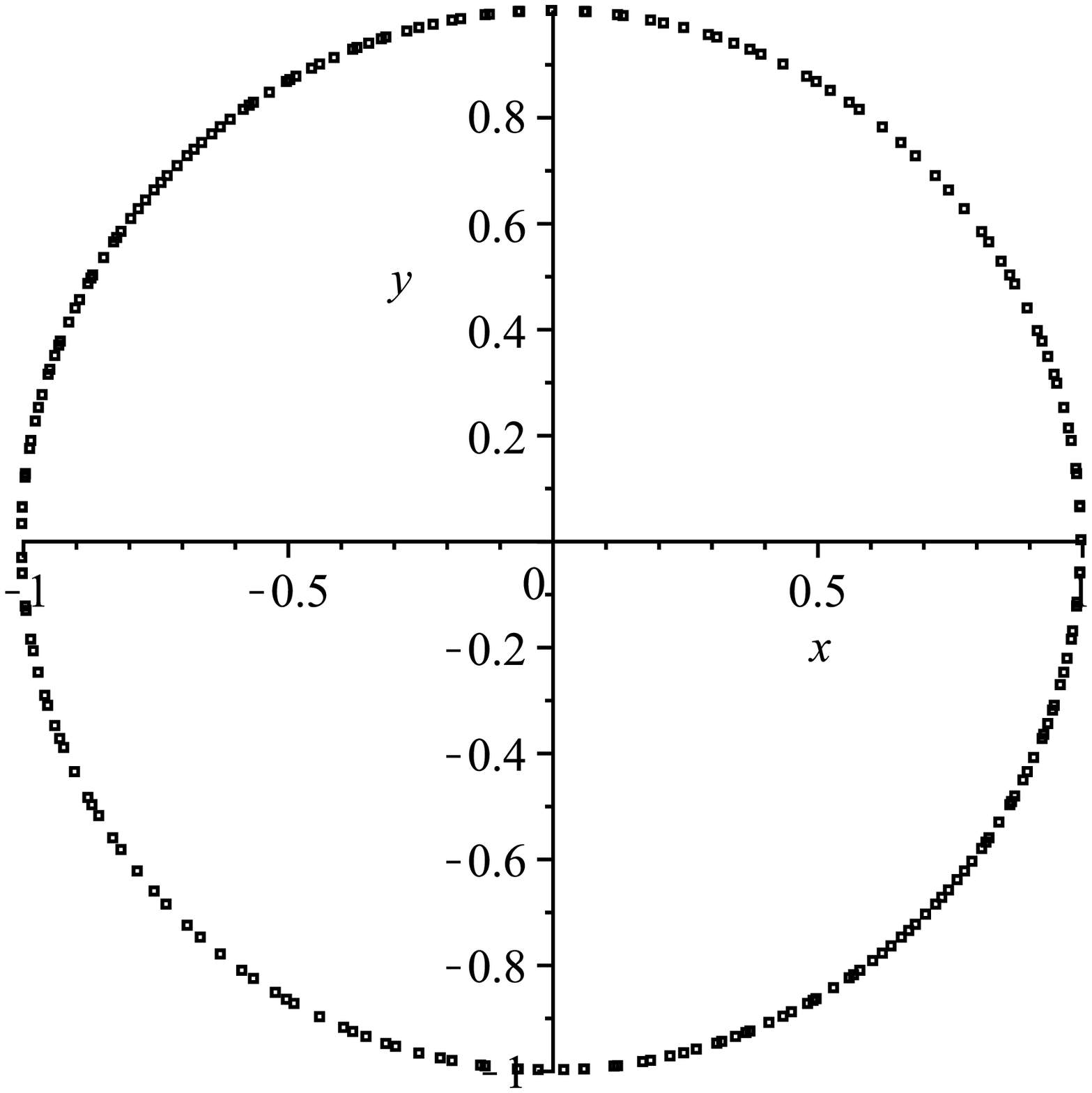,width=5cm,height=5cm}\epsfig{figure=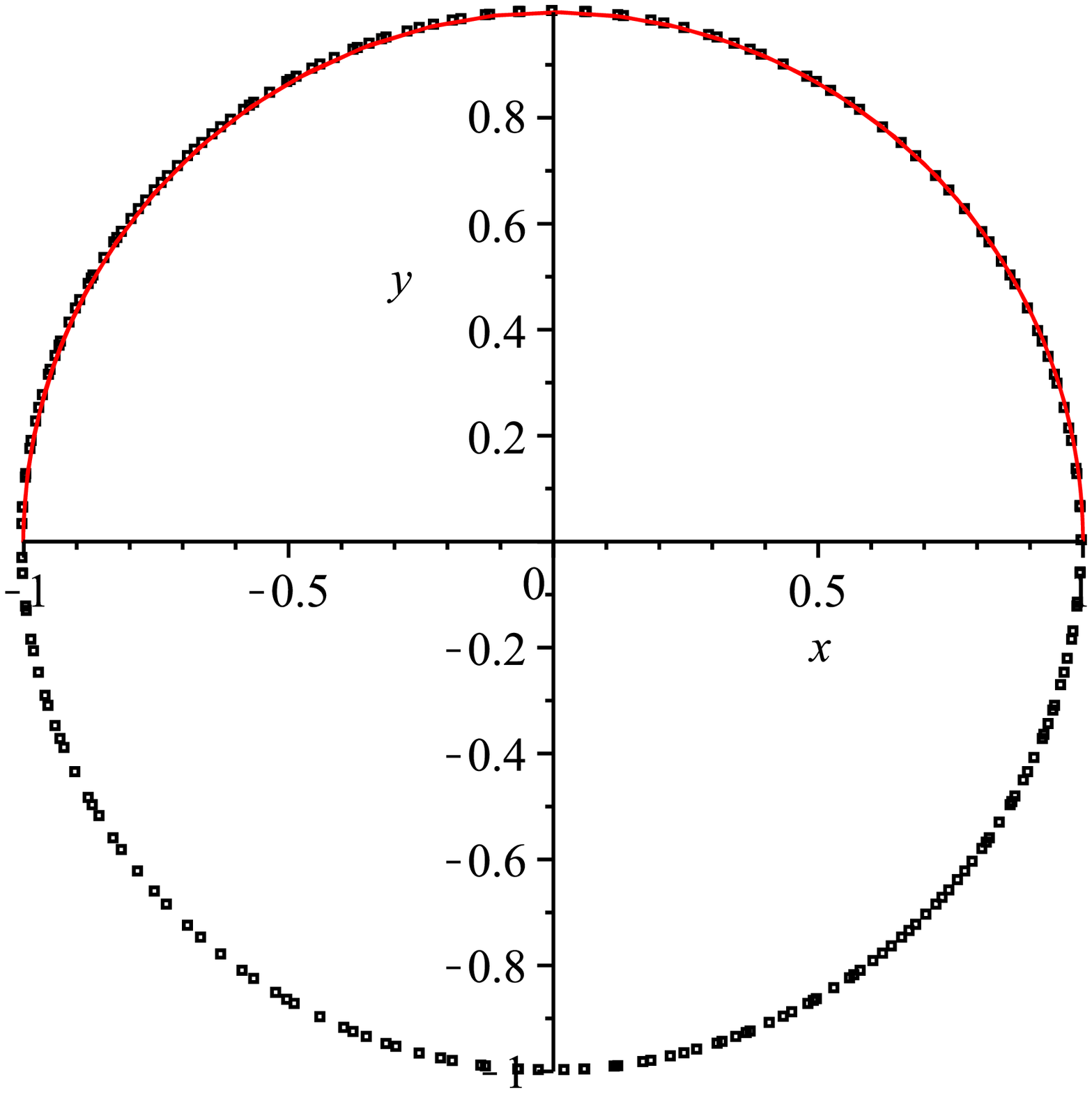,width=5cm,height=5cm}
}}
\end{center}
\vspace*{-1.5cm}
\caption{Input curve $\mathcal C$ (left),  curve ${\mathcal D}$ (center),  curves $\mathcal C$ and ${\mathcal D}$ (right)}
 \label{Ex2fig2}
\end{figure}

\noindent
 One may check that the equality \[L_k(s, x_k)= (x_kq_{k,2}(s)-q_{k,1}(s))^{\ell}+\epsilon^{\ell} W_k(s, x_k),\,\,\|\num(W_k(R,p_k))\|\leq \|{H}^{\cP{\cal Q}}_{k}\|^\ell,\,\,\,k=1,2\] holds. Then, ``{\tt ${\cal
Q}$ is an $\epsilon$-proper reparametrization of $\cal P$}" (see Theorem \ref{Th-resultant} and Corollary \ref{C-resultant}).

\para

 In the following, we analyze the error in the computation, by using Theorem \ref{Th-error}. For this purpose, taking into account the assumptions introduced before Theorem \ref{Th-error}, we consider $I=(-1, 1)$. Thus,  $d=1$. Let  $M\in {\Bbb N}$ be  such that for every $t_0 \in I$, it holds that
$ |{q}_{i, 2}(R(t_0))|\geq M$, and $|p_{i, 2}(t_0)|\geq M$, for $i=1,2$. We have that $M=.2492042100.$ Then, by Theorem \ref{Th-error}, we get that
   $$C= \ell^{1/\ell}{\deg(\cP)}^{1/\ell}=2.828427125,$$ and  for every  $t_0 \in I$, it holds that
   $$|p_i(t_0)-{q}_i(R(t_0))| <  2/M^2 \epsilon   \,C  \|p\|\|q\|=0.9108864449,\quad i=1,2,$$
   where $\|p\|=\|q\|=1$.\\

\para

\noindent
In Figure~\ref{Ex2fig3}, we plot curves $\cal C$ and   $\cal D$ defined by ${\cal P}(t)$ and  ${\cal Q}(t)$, respectively, for $t\in I$.

\begin{figure}[h]
\begin{center}
\hbox{
\centerline{\epsfig{figure=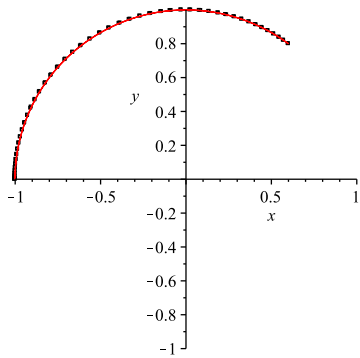,width=5cm}
}}
\end{center}
\caption{Curves $\mathcal C$ and ${\mathcal D}$ for $t\in I$}
\label{Ex2fig3}
\end{figure}
 \end{example}

\para

\begin{example}   Let   $\epsilon=0.0001$, and the rational curve $\cal C$  defined by the parametrization \[{\cal P}(t)=\left(\frac{p_{1,1}(t)}{p_{1,2}(t)},\,
\frac{p_{2,1}(t)}{p_{2,2}(t)}\right)=\]\[\left(\frac{.7498125469 t^6+t^3+.4973756561}{1.749562609t^6+1.749812547t^3+.2499375156}, \frac{.0002499375156 t (10000 t^5+1.)}{17.49562609 t^6+17.49812547 t^3+2.499375156}\right).\]
Using the
{\it SNAP}   package, one has that  $\egcd(p_{j, 1}, p_{j, 2})=1,\,\,j=1,2$. We apply the algorithm and, in Step 1, we compute the polynomials\\

\noindent $H_{1}^{\cP\cP}(t, s)=-7004000 s^6 t^3-10930000 s^6+7004000 s^3 t^6-9930990 s^3+10930000 t^6+9930990 t^3,$\\

\noindent $H_{2}^{\cP\cP}(t,s)= 70010000 s^6 t^3+10000000 s^6+7000 s t^6+7001 s t^3+1000s-70010000 s^3 t^6-10000000 t^6-7000 t s^6-7001ts^3-1000t.$\\

\noindent Now, we compute the polynomial   $S^{{\cal P}{\cal P}}_\epsilon$. We have,
$$S^{{\cal P}{\cal P}}_\epsilon(t, s)\approx_\epsilon C_0(t)+C_1(t)s+C_2(t)s^2+C_3(t)s^3,$$
where $$C_0(t)=t (69970939184+535492598272100802900 t^2)
,\quad $$$$C_1(t)= t (-52478204388+535492598272100802900 t)-69970939184-535492598272100802900 t^2
,\quad $$$$C_2(t)=63943722313 t+52478204388,\quad C_3=-535492598336044525213.$$
 Then,  $\eindex({\cal P})=\deg_t(S^{{\cal P}{\cal P}}_\epsilon)=2$ (see Definition \ref{Def-approxindex}). Now, we apply Step 4 of the algorithm, and we consider
\[R(t)= \frac{C_0(t)}{C_3(t)}=\frac{-4t (17492734796+133873149568025200725 t^2)}{535492598272100802900}.\]
In Steps 5 and 6 of the algorithm, we determine the polynomials $L_k(s,x_k)$, and  we compute the root in the variable $x_k$ of the polynomial $\frac{\partial L_k}{\partial x_k}(s, x_k),\,k=1,2$. We get the rational parametrization  $\widetilde{\cal Q}(t)$. We simplify it, and we return the curve $\cal D$ defined by the $\epsilon$-numerical reparametrization
  $$ {\cal Q}(t)=\left(\frac{.7498125351 t^2+t+.4973756559}{1.749562581 t^2+1.749812559 t+.2499375114},\right.$$$$\left. \frac{.2499375117 t^2+.5551941368\,10^{-8}t-.5495954487\,10^{-8}}{1.749562581 t^2+1.749812559 t+.2499375114}\right).$$

  \para

 One may check that the equality of Theorem \ref{Th-resultant} does not hold. However, Remark \ref{R-Thresultant} holds taking $\overline{\epsilon}=0.0005$. Then, ``{\tt ${\cal
Q}$ is an $\overline{\epsilon}$-proper reparametrization of $\cal P$}". In Figure~\ref{Ex3fig1}, we plot the input curve $\cal C$ and the output curve $\cal D$.

\begin{figure}[h]
\begin{center}
\hbox{
\centerline{\epsfig{figure=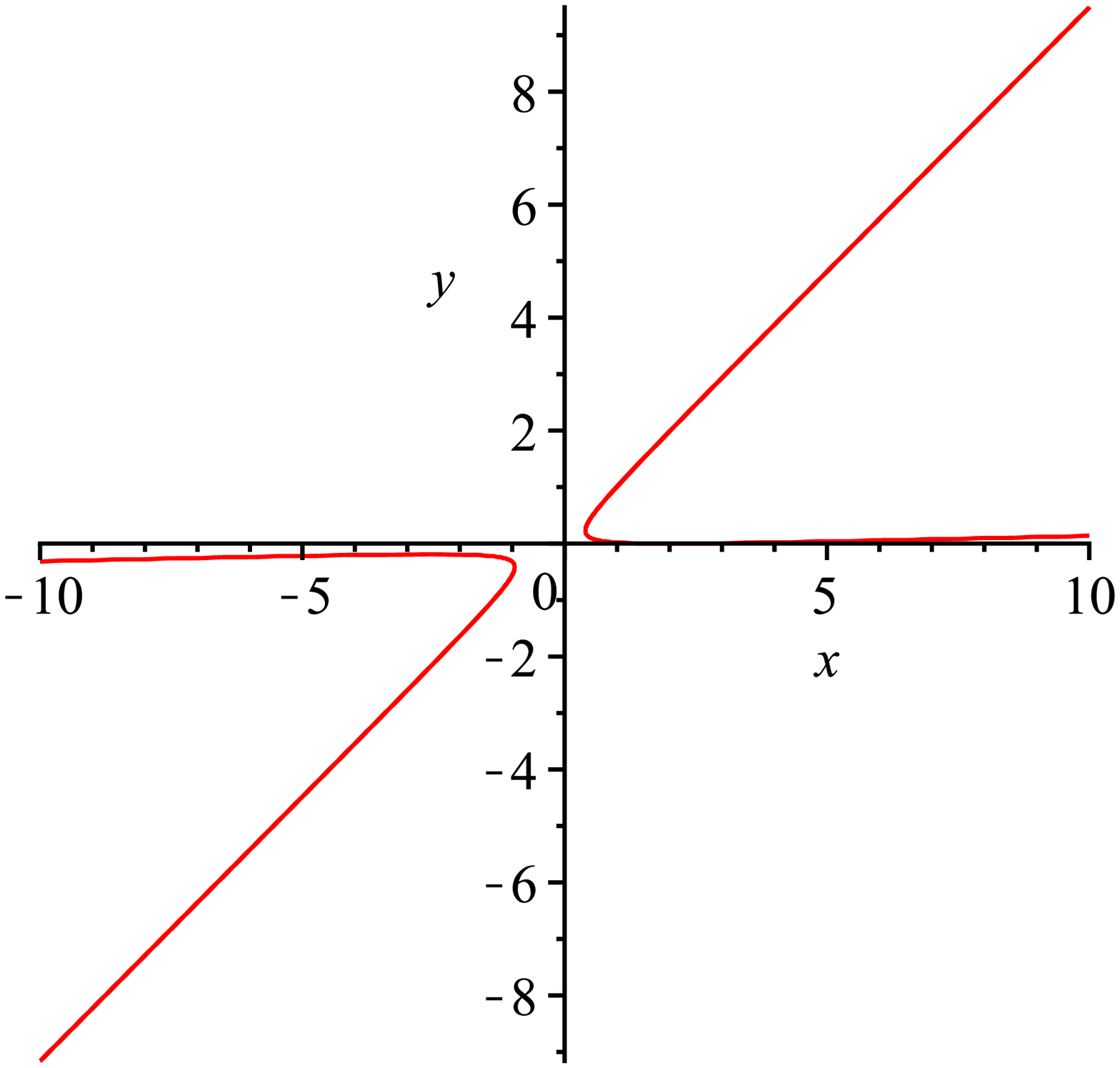,width=5cm,height=5cm} \epsfig{figure=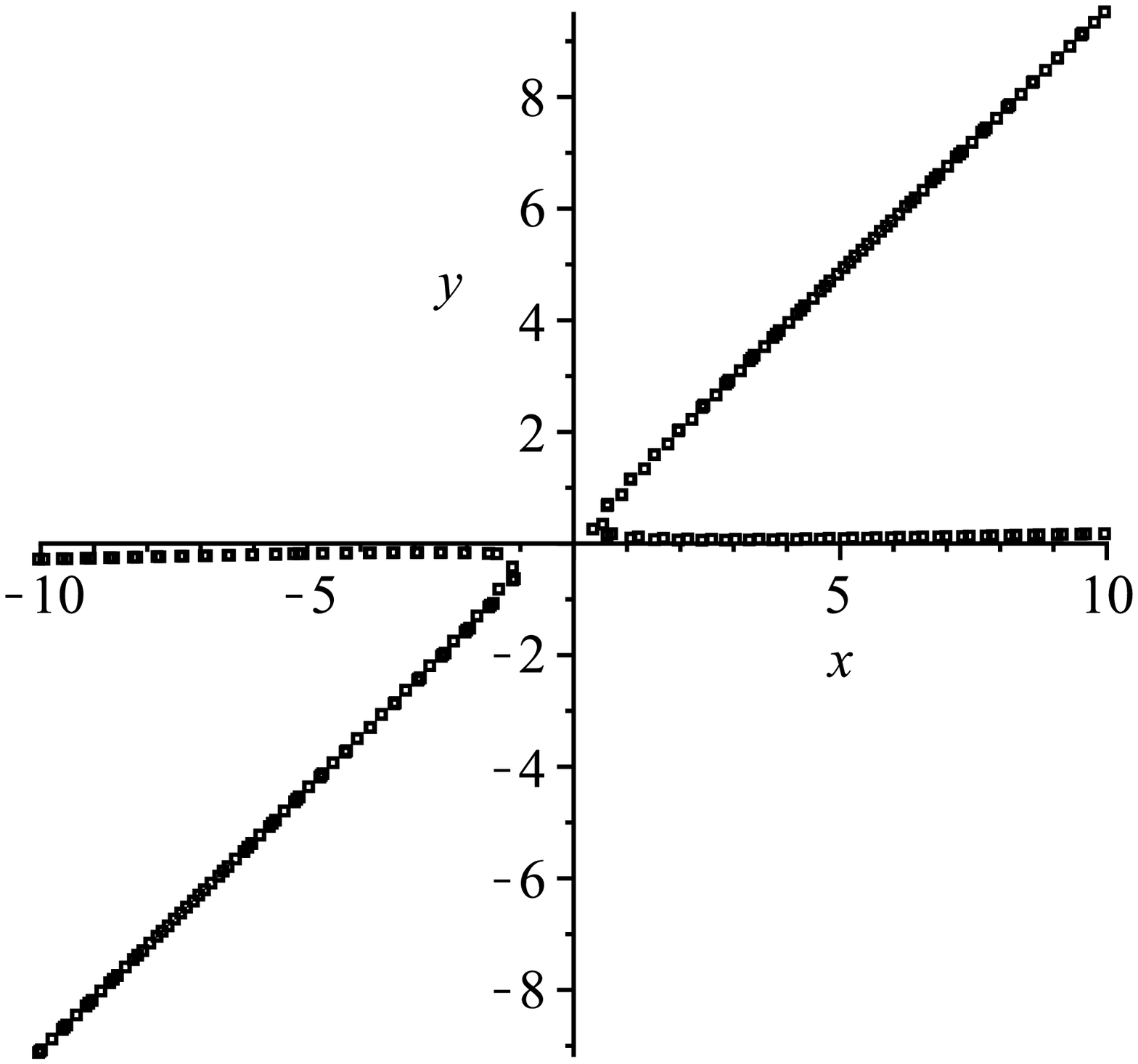,width=5cm,height=5cm} \epsfig{figure=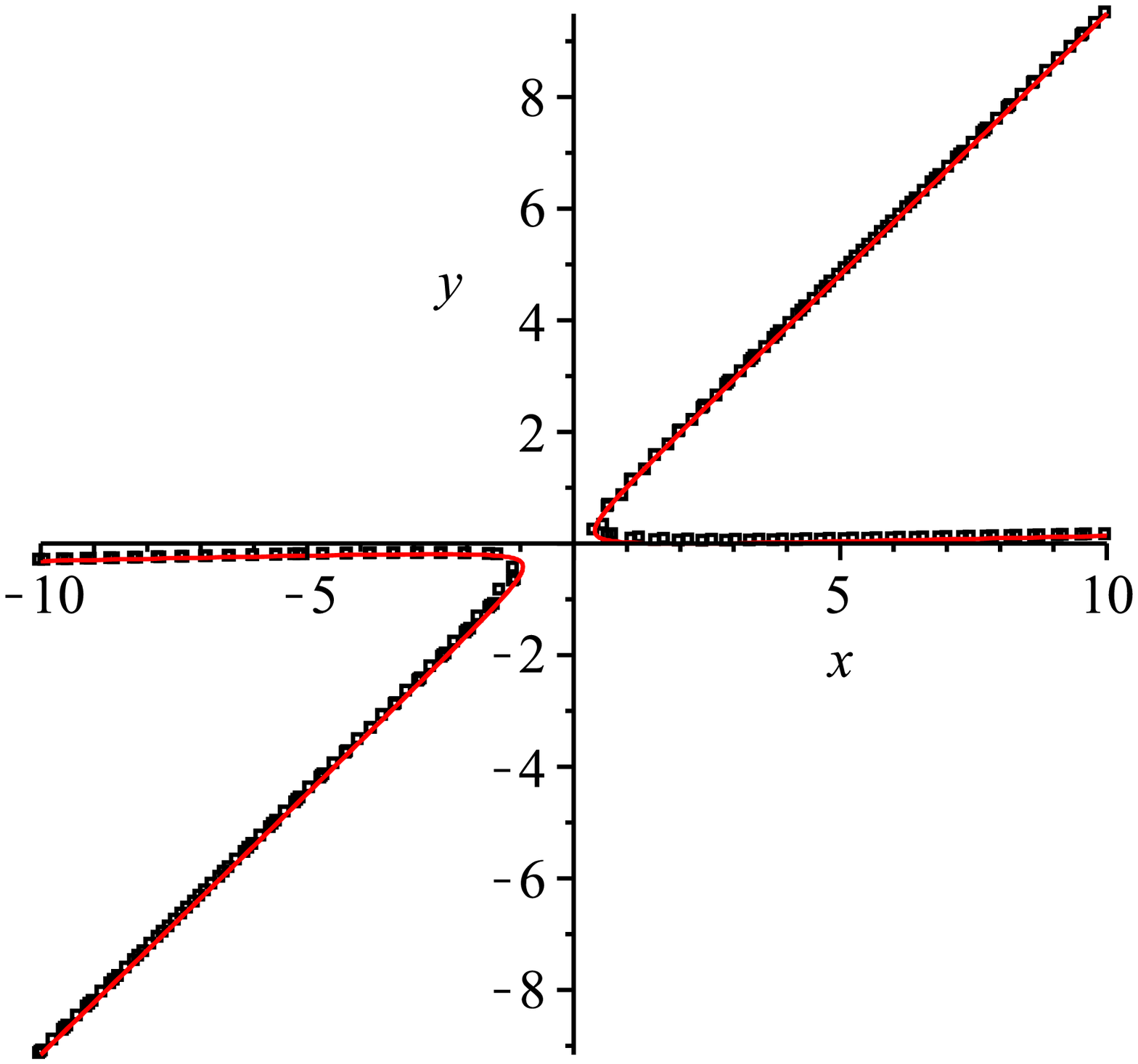,width=5cm,height=5cm}
}}
\end{center}
\caption{Input curve $\mathcal C$ (left),  curve ${\mathcal D}$ (center),  curves $\mathcal C$ and ${\mathcal D}$ (right)}
 \label{Ex3fig1}
\end{figure}

\para
\para

 In the following, we analyze the error in the computation, by using Theorem \ref{Th-error}. For this purpose, we consider $I=(3, 10)$. Thus,  $d=10$. Let  $M\in {\Bbb N}$ be  such that for every $t_0 \in I$, it holds that
$ |{q}_{i, 2}(R(t_0))|\geq M$, and $|p_{i, 2}(t_0)|\geq M$, for $i=1,2$. We have that $M=1322.925998.$ Then, by Theorem \ref{Th-error}, we deduce that
   $$C= \displaystyle{\frac{d^{\deg(\cP)+1}}{(d-1)^{1/\ell}}}=4807498.567,$$ and  for every  $t_0 \in I$, it holds that
   $$|p_i(t_0)-{q}_i(R(t_0))| < 2/M^2 \overline{\epsilon}    \,C\|p\|\|q\|=.08410680133,\quad i=1,2,$$
    where $\|p\|=17.49812547$, and $\|q\|=1.749812559$.\\

\para

\noindent
In Figure~\ref{Ex3fig2}, we plot curves $\cal C$ and   $\cal D$ defined by ${\cal P}(t)$ and  ${\cal Q}(t)$, respectively, for $t\in I$.

\begin{figure}[h]
\begin{center}
\hbox{
\centerline{\epsfig{figure=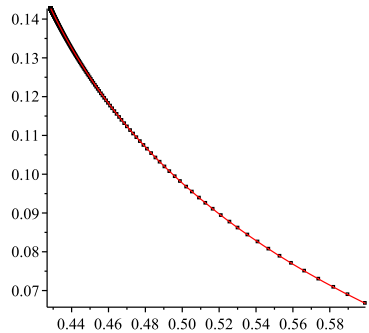,width=5cm}
}}
\end{center}
\vspace*{-1.5cm}
\caption{Curves $\mathcal C$ and ${\mathcal D}$ for $t\in I$}
\label{Ex3fig2}
\end{figure}
 \end{example}

 \para

\begin{example}
\begin{center}
\begin{tabular}{|lll|}
\hline
Input Curve $\mathcal C$ & \vline & ${\cal P}=({\frac {{t}^{6}-3{t}^{5}- 3.001{t}^{4}+ 11.001{t}^{3}+9{t}^{2
}-15t- 9.002}{{t}^{2}-t- 2.001}},{\frac {{t}^{4}- 2.001{t}^{3}-2
{t}^{2}+ 3.002t+3}{{t}^{2}-t- 2.001}}
)$\\
\hline
Tolerance $\epsilon$ & \vline & $0.02$ \\
\hline
Output Curve $\cal D$ & \vline &
${\cal Q}=({\frac { 0.06667333664\,{t}^{3}- 0.4000900188\,{t}^{2}+t-
 0.6001100173}{ 0.06667333664\,t- 0.1334077982}},$\\
 & \vline &
 \hspace*{1.0cm} ${\frac {
 0.06667333662\,{t}^{2}- 0.2001089149\,t+ 0.2000366790}{ 0.06667333664
\,t- 0.1334077982}}
)$ \\
\hline
Error bound & \vline &
$|p_i(t_0)-{q}_i(R(t_0))| \leq 0.4582153762,$  $ \qquad I=(0, 0.5)$  \\
\hline
Figures & & \\
\epsfig{figure=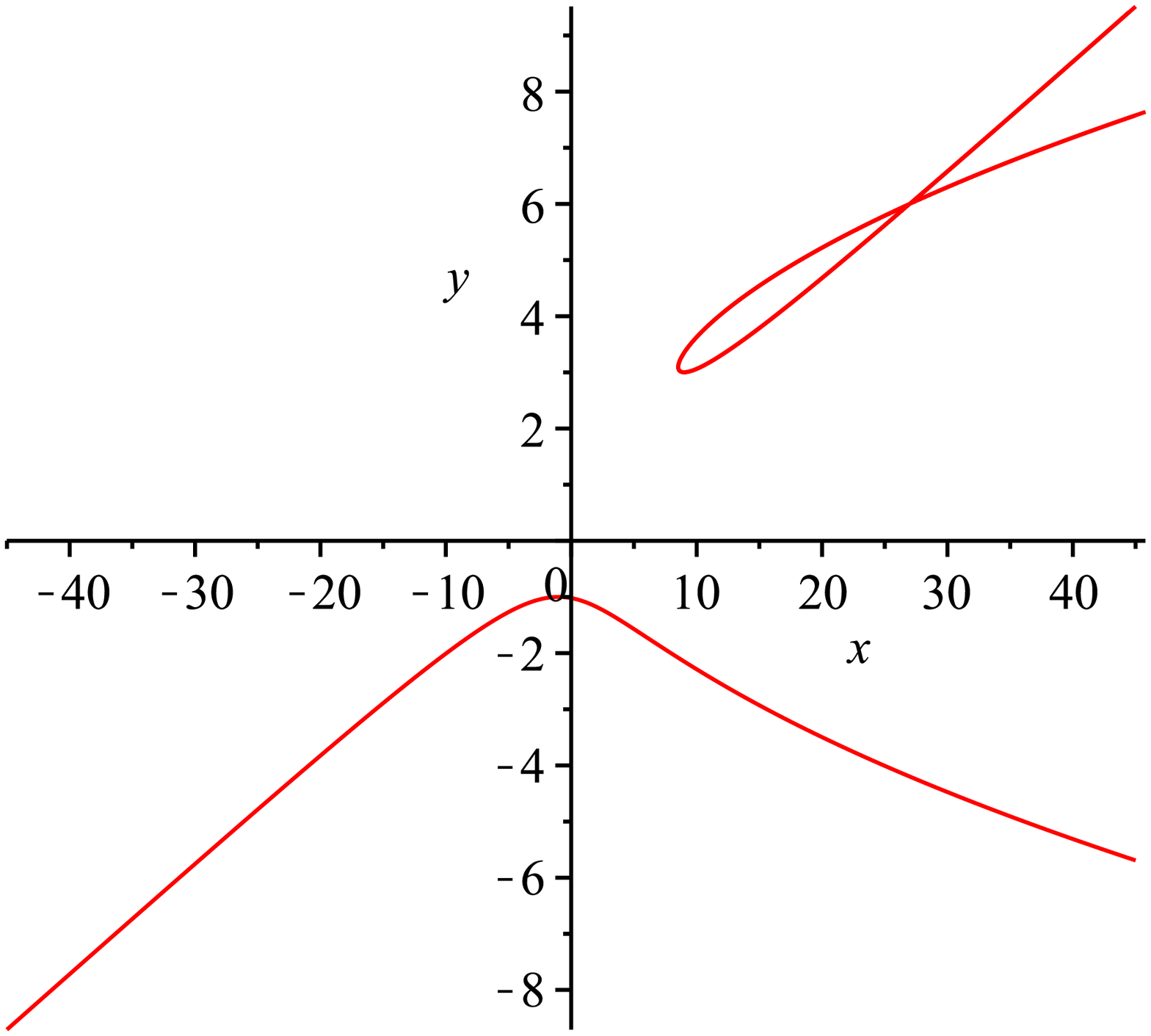,width=4.5cm} & & \hspace*{0.2cm}
\epsfig{figure=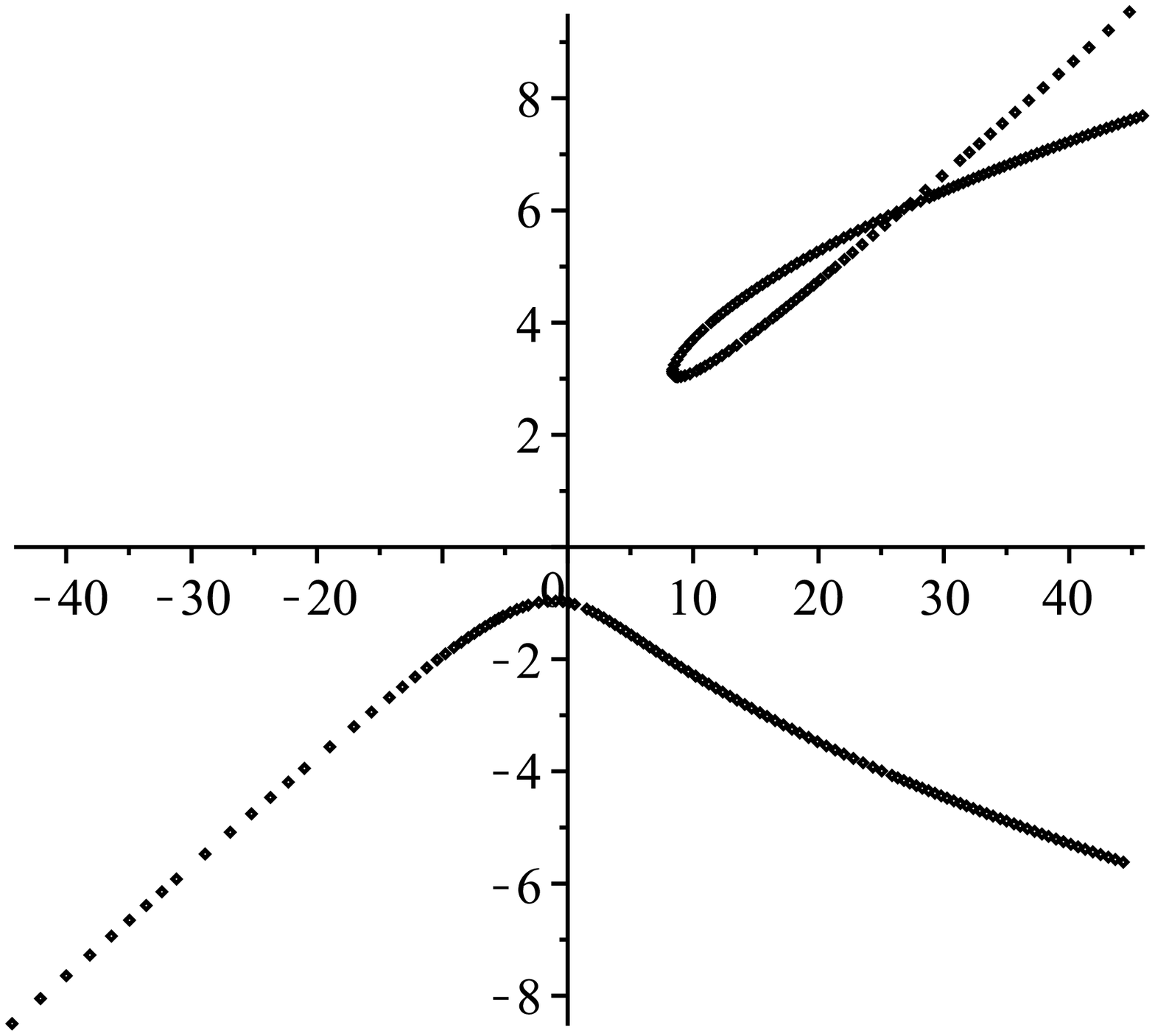,width=4.5cm}  \hspace*{1.0cm}
\epsfig{figure=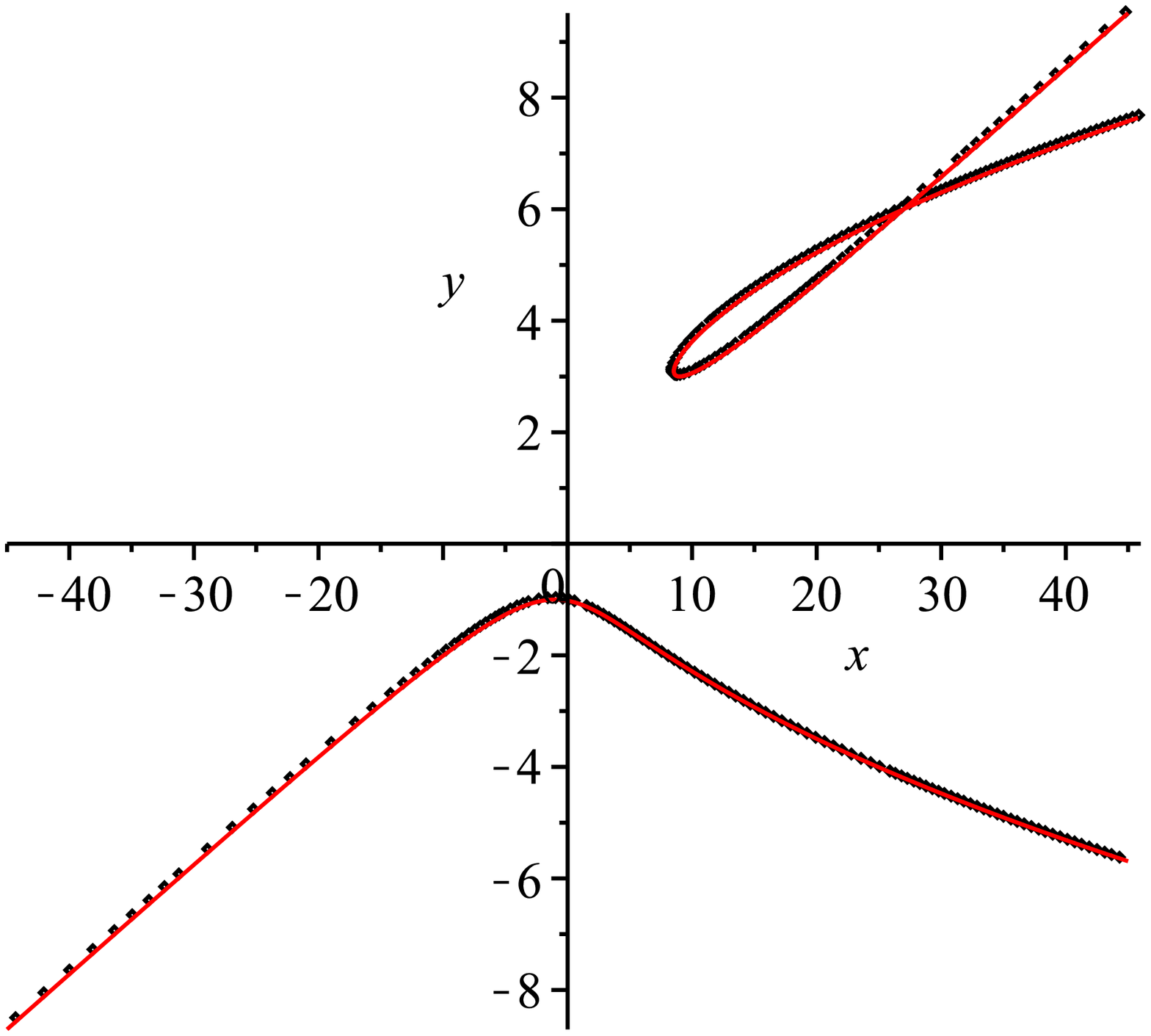,width=4.5cm}
\\
\hspace*{1.1cm} \mbox{Curve}\,\, $\mathcal C$& &\hspace*{1.1cm} \mbox{Curve}\,\,\, {$\mathcal D$}\hspace*{3.1cm}  \mbox{Curves}\,\, $\mathcal C$ \mbox{and} $\mathcal D$\\
\hline
\end{tabular}
\end{center}

 \end{example}

 \para

\begin{example}
\begin{center}
\begin{tabular}{|lll|}
\hline
Input Curve $\mathcal C$ & \vline & ${\cal P}=({\frac {.4995004995 t^6-t^3+.5005024975}{.02497502498t^8-.2248001998t^6+.4495504496t^3-.2247752248}}
,$\\
 & \vline &
 \hspace*{0.2cm} $
{\frac {1000t^{12}-0.2}{(1000t^8-9001t^6+18000t^3-9000)(t^3-1)}}
)$\\
\hline
Tolerance $\epsilon$ & \vline & $0.01$ \\
\hline
Output Curve $\cal D$ & \vline &
${\cal Q}=(\frac{.4490987986\,10^{-11}t^2-.1511966995\,10^{-5}t+1}{.04989990076t^2+.3006866597\,10^{-6}t-.4112066119}
,$\\
 & \vline &
 \hspace*{0.2cm} $\frac {.04989241582t^3+.2507815761\,10^{-6}t^2+.03788675995t+.4671769265\,10^{-5}}{{.04989990076t^2+.3006866597\,10^{-6}t-.4112066119}}
)$ \\
\hline
Error bound & \vline &
$|p_i(t_0)-{q}_i(R(t_0))| \leq 0.6913866372,$  $ \qquad I=(-2,2)$  \\
\hline
Figures & & \\
\epsfig{figure=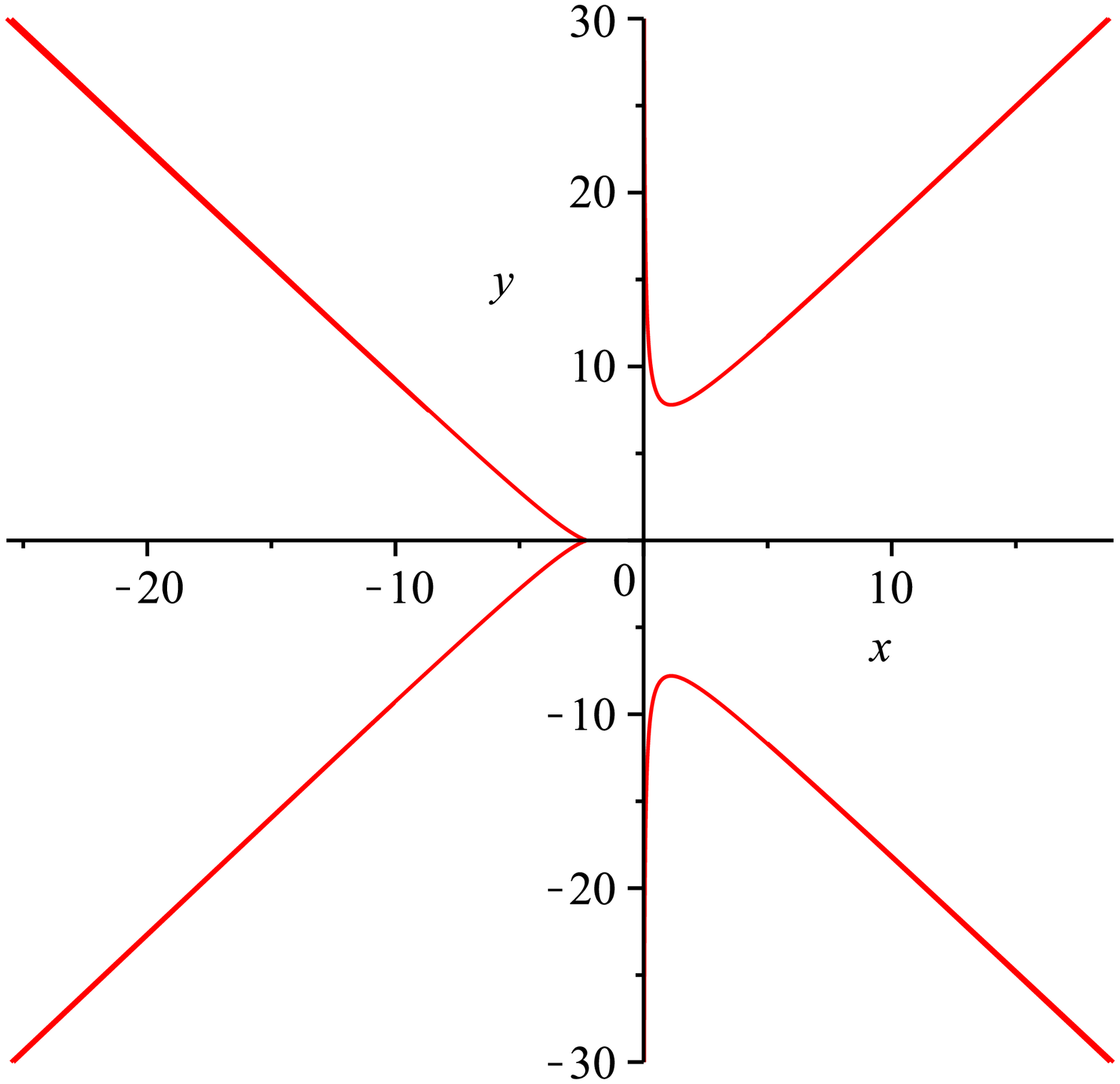,width=4.5cm} & & \hspace*{0.2cm}
\epsfig{figure=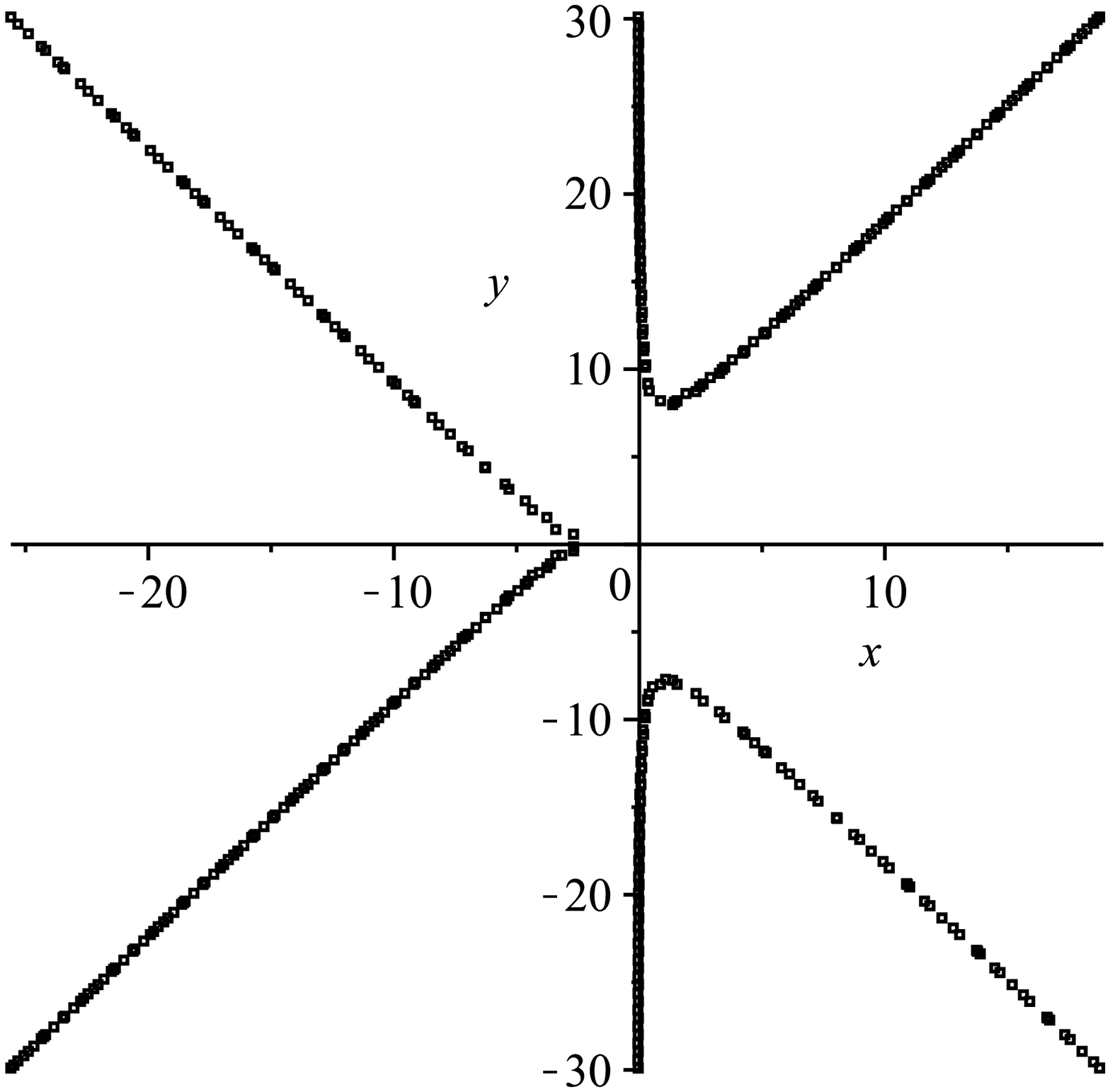,width=4.5cm}  \hspace*{1.0cm}
\epsfig{figure=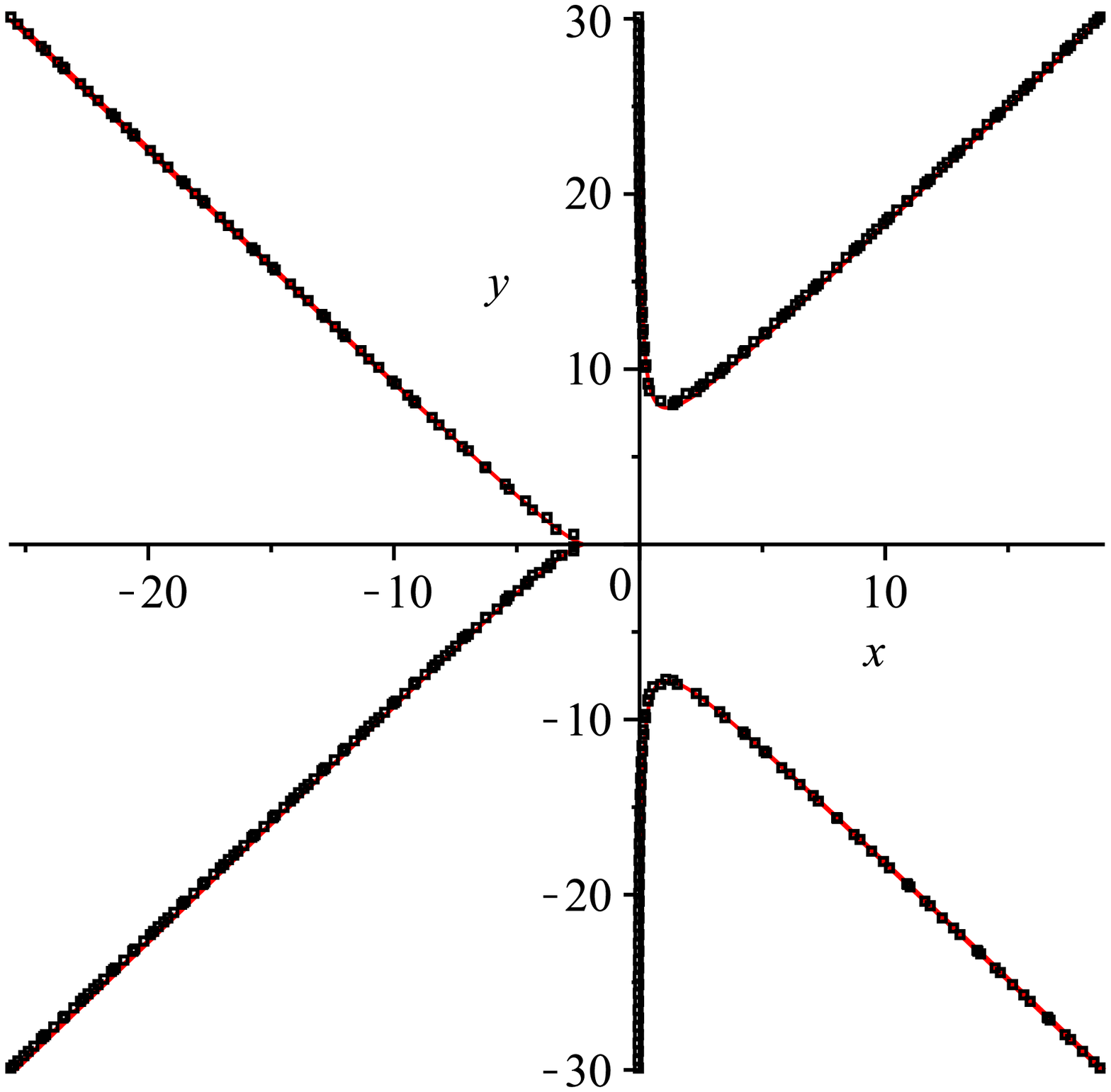,width=4.5cm}
\\
\hspace*{1.1cm} \mbox{Curve}\,\, $\mathcal C$& &\hspace*{1.1cm} \mbox{Curve}\,\,\, {$\mathcal D$}\hspace*{3.1cm}  \mbox{Curves}\,\, $\mathcal C$ \mbox{and} $\mathcal D$\\
\hline
\end{tabular}
\end{center}

 \end{example}

\newpage

\begin{example}
\begin{center}
\begin{tabular}{|lll|}
\hline
Input Curve $\mathcal C$ & \vline & ${\cal P}=({\frac { 20.001\,{t}^{8}-40\,{t}^{5}+20\,{t}^{2}+2\,{t}^{7}- 2.001\,{
t}^{4}-{t}^{6}}{ \left( {t}^{3}- 1.001 \right) ^{3}}}$,
\\
 & \vline & \hspace*{1.0cm} ${\frac {-2\,{t}^
{4}- 6.002\,{t}^{5}+6\,{t}^{2}+6\,{t}^{6}- 12.002\,{t}^{3}+ 6.002}{
 \left( -{t}^{2}+{t}^{3}-1 \right)  \left( {t}^{3}- 1.001 \right) }}
)$\\
\hline
Tolerance $\epsilon$ & \vline & $0.001$ \\
\hline
Output Curve $\cal D$ & \vline &
${\cal Q}=({\frac {-{t}^{2}+ 0.1002417588\,t+ 0.04999508896}{
 0.050095827\,{t}^{3}- 1.520248235\,10^{-5}\,{t}^{2}+
 4.9502838\,10^{-8}\,t- 2.19451487\,10^{-10}}},
$\\
 & \vline &
 \hspace*{1.0cm} ${\frac {9.046219880\,10^{-4}\,{t}^{2}+ 0.0009044005499\,t- 0.0003016044631}{
 0.0001508373032\,{t}^{2}+ 0.0001508171174\,t- 1.509849726\,10^{-7}}}
)$ \\
\hline
Error bound & \vline &
$|p_i(t_0)-{q}_i(R(t_0))| \leq0.1254659264,$  $ \qquad I=(-5,5)$  \\
\hline
Figures & & \\
\epsfig{figure=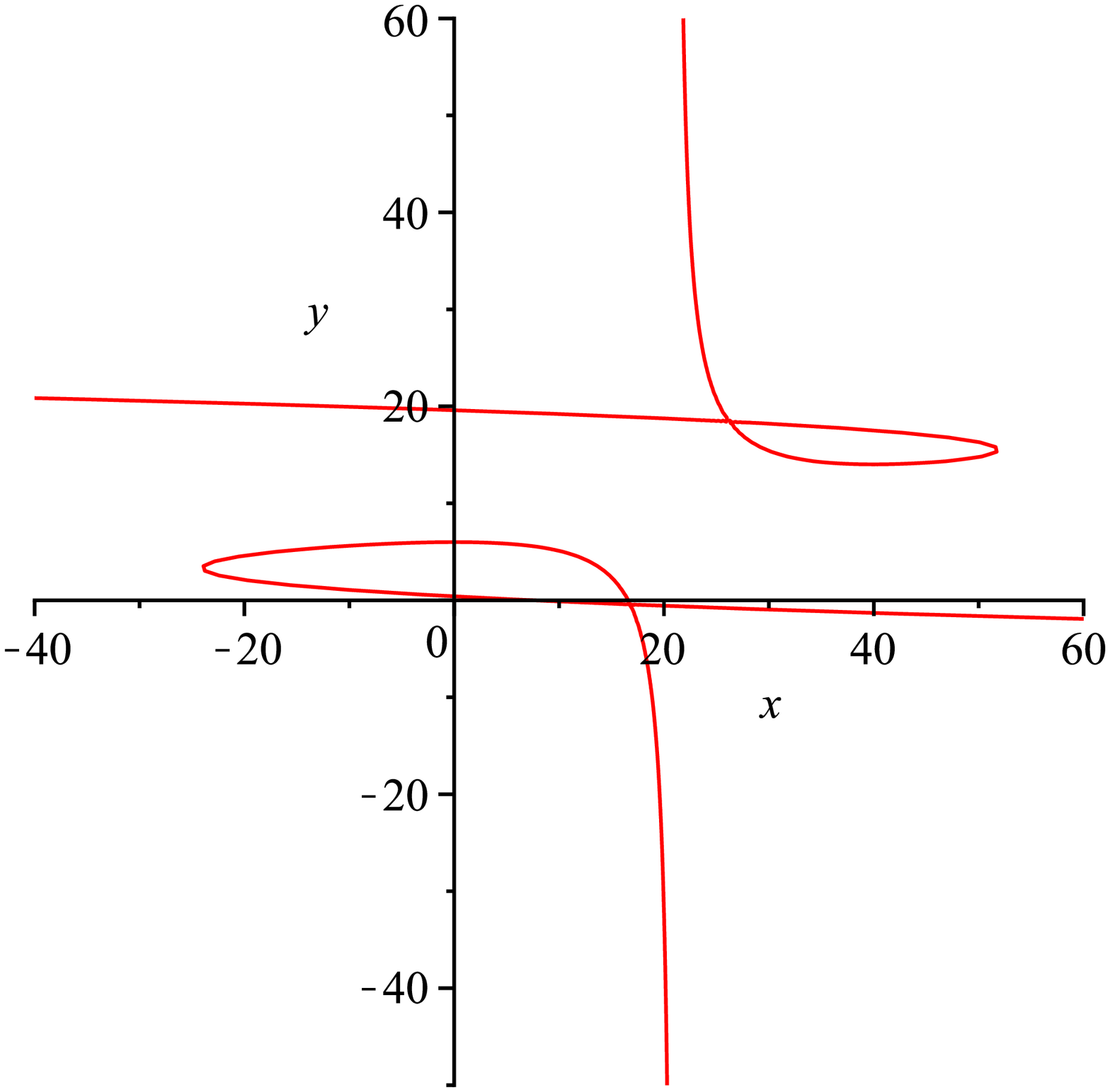,width=4.5cm} & & \hspace*{0.2cm}
\epsfig{figure=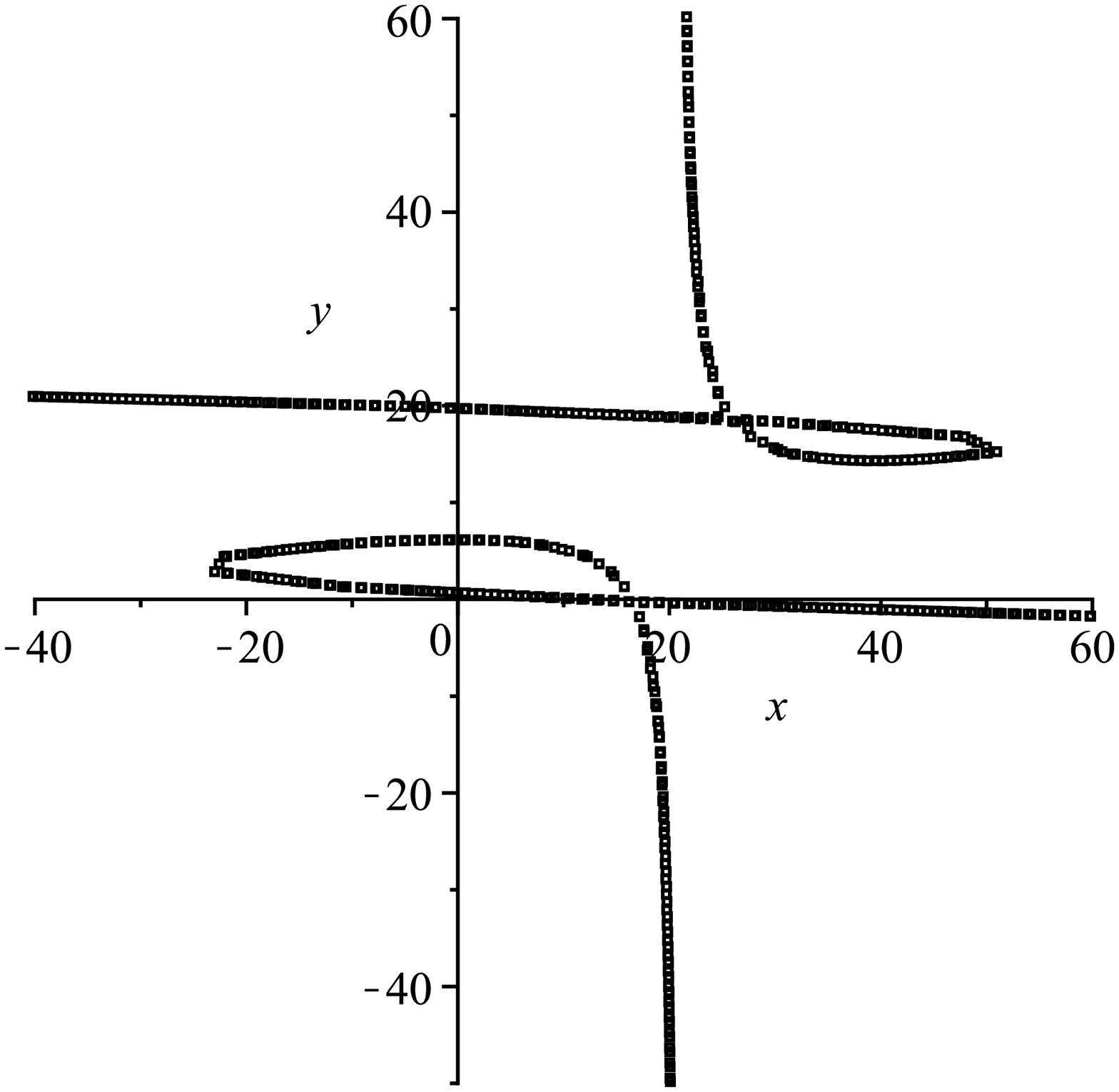,width=4.5cm}  \hspace*{1.0cm}
\epsfig{figure=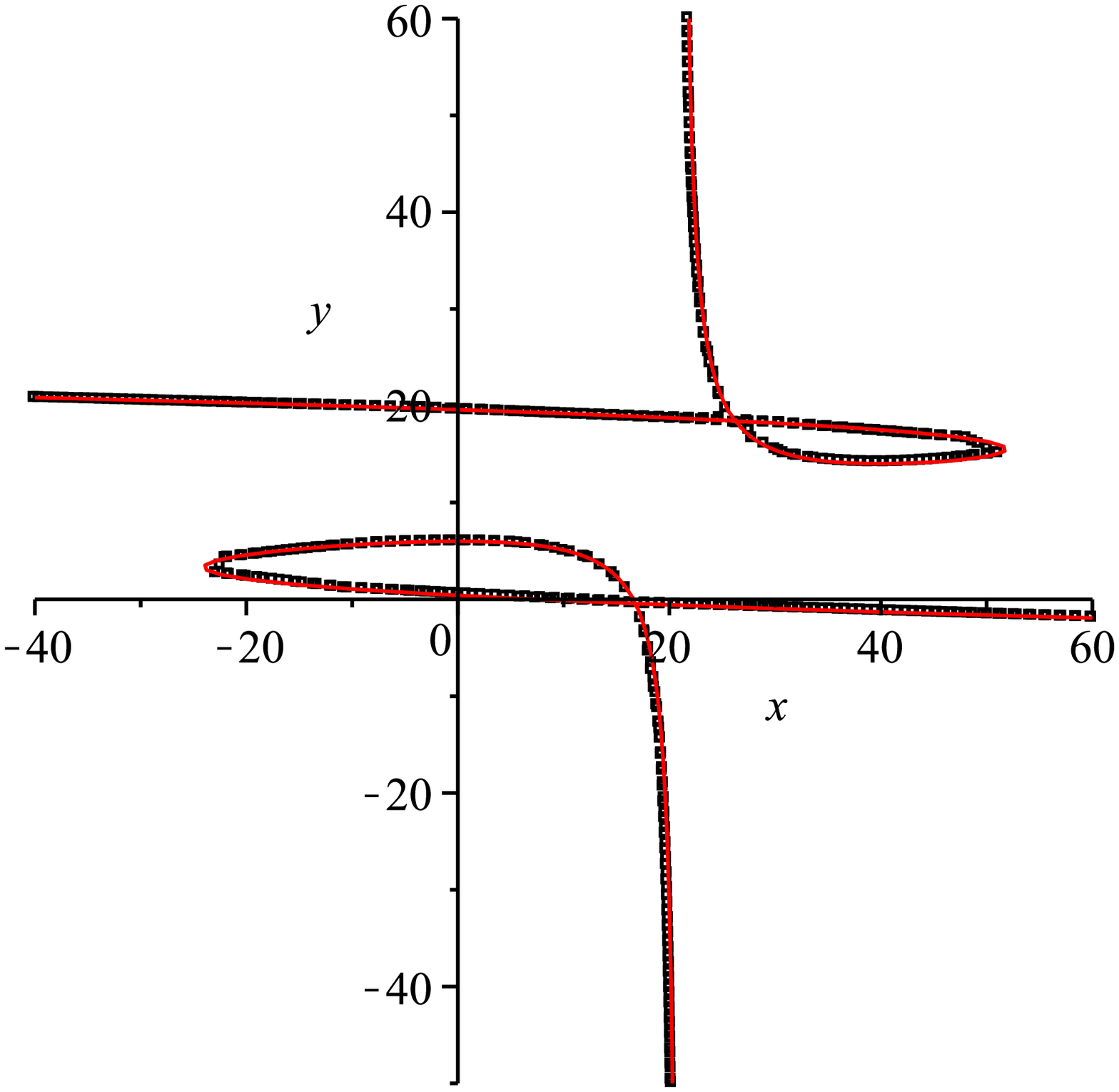,width=4.5cm}
\\
\hspace*{1.1cm} \mbox{Curve}\,\, $\mathcal C$& &\hspace*{1.1cm} \mbox{Curve}\,\,\, {$\mathcal D$}\hspace*{3.1cm}  \mbox{Curves}\,\, $\mathcal C$ \mbox{and} $\mathcal D$\\
\hline
\end{tabular}
\end{center}

 \end{example}


\section{Conclusion}
The paper focus on the problem of numerical proper reparametrization which has both theoretical and practical background. Based on the results on the symbolic situation (see \cite{Perez-repara}), we achieve the expected properties and algorithm. For a given numerical curve, we can determine whether it is approximate improper with respect to a given precision  and, in the affirmative case, an $\epsilon$-proper reparametrization can be found. More important, the reparameterized curve obtained always lies in the certain offset region of the input one (and reciprocally). As a natural but more difficult problem, we would like to consider the problem of numerical proper reparametrization for rational space curves and surfaces.\\



\begin{thebibliography}{20}

\bibitem{Abh} Abhyankar, S., Bajaj, C.  (1988). {\it Automatic Parametrization of Rational curves and Surfaces III: Algebraic Plane Curves}. Computer Aided Geometric Design. Vol.  5, pp. 321--390.

\bibitem{ASS} Arrondo, E., Sendra, J.,
Sendra, J.R.  (1997). {\it Parametric Generalized Offsets to
Hypersurfaces}. J. Symbolic Computation. Vol.   23, pp. 267-285.


\bibitem{Beckermann1} Beckermann, B., Labahn, G.  (1998). {\it A fast and numerically stable Euclidean-like algorithm for detecting relatively prime numerical polynomials.} Journal of Symbolic Computation. Vol. 26, pp. 691-714.

\bibitem{Beckermann2} Beckermann, B.,  Labahn, G. (1998). {\it When are two numerical polynomials relatively prime?}  Journal of Symbolic Computation. Vol. 26, pp. 677-689.



\bibitem{chionh06}
 Chionh, E.-W.,  Gao, X.-S.,   Shen,  L.-Y. (2006). {\it Inherently improper surface parametric supports}.  Computer Aided Geometric Design. Vol.  23(8), pp. 629-639.


    \bibitem{Cor2} Corless, R.M., Giesbrecht, M.W., Kotsireas, I.S. van Hoeij, M., Watt, S.M. (2001).
 {\it Towards Factoring Bivariate Approximate Polynomials}. Proc.
 ISSAC 2001. pp. 85--92.


\bibitem{Corless} Corless, R.M., Watt, S. M., Zhi, L. (2004). \emph{QR Factoring to compute the GCD of Univariate Approximate Polynomials.} IEEE Transactions on Signal Processing. Vol. 52(12), pp.  3394--3402.


\bibitem{Cox1998}  Cox, D.A.,
Little, J.,  O'Shea, D.  (1998). {\it Using Algebraic Geometry. Graduate
Texts in Mathematics}. Vol.  185. Springer--Verlag.

\bibitem{CLO2} Cox, D.A., Sederberg, T.W., Chen
F., (1998). {\it The moving line ideal basis of planar rational
curves.} Computer Aided Geometric Design. Vol.  8, pp. 803--827.

\bibitem{Farouki} Farouki, R.T.,   Rajan, V.T. (1988). {\it On the numerical condition of algebraic curves and surfaces. 1: Implicit
equations}. Computer Aided Geometric Design. Vol.  5, pp. 215-252.

\bibitem{galliao02}Galligo, A., Rupprech, D. (2002). \emph{Irreducible decomposition of curves}. Journal of Symbolic Computation. Vol.33, pp. 661¨C677.

\bibitem{gao92}
 Gao, X.-S., Chou, S.-C. (1992). {\it Implicitization of rational parametric equations}.  Journal of Symbolic Computation. Vol.14(5), pp. 459--470.

\bibitem{HSW} Hoffmann, C.M., Sendra, J.R., Winkler, F. (1997).
{\it Parametric Algebraic Curves and Applications}. J. Symbolic
Computation. Vol. 23.

\bibitem{HL97} Hoschek, J., Lasser, D. (1993). {\it Fundamentals
of Computer Aided Geometric Design}.  A.K. Peters Wellesley MA.,
Ltd.

\bibitem{erich08}Kaltofen, E., May, J.P., Yang, Z., Zhi, L. (2008). \emph{ Approximate factorization of multivariate polynomials using singular value decomposition}, Journal of Symbolic Computation, Vol.(43)5, pp. 359-376.

\bibitem{Karmarkar} Karmarkar, N.,  Lakshman, Y.N. (1996).  {\it Approximate polynomial greatest common divisors and nearest singular polynomials.} ISSAC 1996, pp. 35-39. ACM Press.


\bibitem{Perez-repara} P\'erez-D\'{\i}az, S. (2006).
{\it On the Problem of Proper Reparametrization for Rational Curves and Surfaces}. Computer Aided Geometric Design. Vol. 23(4), pp.  307-323.


\bibitem{diaz02} P\'erez-D\'{\i}az, S., Schicho, J., Sendra, J.R. (2002).
{\it Properness and inversion of rational parametrizations of
surfaces}.  Appl. Algebra Eng. Commun. Comput. Vol.  13(1), pp. 29--51.


 \bibitem{PSS} P\'erez-D\'{\i}az, S., Sendra, J.R., Sendra, J. (2004).
 {\it  Parametrizations of Approximate Algebraic Curves by Lines}. Theoretical Computer Science on  Algebraic - Numeric  Algorithms. Vol.  315/2-3, pp. 627-650.


\bibitem{PSS1} P\'erez-D\'{\i}az, S., Sendra, J.R., Sendra, J. (2005).
 {\it  Parametrizations of Approximate Algebraic Surfaces by Lines}. Computer Aided Geometric Design. Vol. 22(2). pp. 147-181.


\bibitem{PSS2} P\'erez-D\'{\i}az, S., Sendra, J.R., Sendra, J. (2006).
 {\it  Distance Bounds of e-Points on Hypersurfaces}. Theoretical Computer Science. Vol 359. N. 1-3. pp. 344-368.


\bibitem{PSV} P\'erez-D\'{\i}az, S., Sendra, J.R., Villarino, C. (2007).
{\it Finite Piecewise Polynomial Parametrization of Plane Rational Algebraic Curves}. Applicable Algebra in Engineering, Communication and Computing. Vol 18(1-2), pp.  91-105.

\bibitem{PSS3} P\'erez-D\'{\i}az, S., Rueda, S.L., Sendra, J.R., Sendra J. (2010).
 {\it  Approximate parametrization of plane algebraic curves by linear systems of curves}. Computer Aided Geometric Design. Vol 27. pp. 212-31.



\bibitem{Sed86} Sederberg, T.W. (1986). {\it Improperly
Parametrized Rational Curves}. Computer Aided Geometric Design. Vol.  3, pp.  67-75.

\bibitem{Sen2} Sendra, J.R., Winkler, F. (2001). {\it Tracing Index of Rational Curve Parametrizations.}
Computer Aided Geometric Design. Vol.  18(8), pp.  771-795.


\bibitem{SWP}  Sendra, J.R., Winkler, F., Perez-Diaz, S. (2007). {\it Rational Algebraic Curves: A Computer Algebra Approach}. Series: Algorithms and Computation in Mathematics.  Vol. 22. Springer Verlag.


\bibitem{shen06}
 Shen, L.-Y.,  Chionh, E.-W.,  Gao, X.-S.,   Li, J. (2011).
{\it Proper reparametrization for inherently improper unirational varieties}. Journal of Systems Sciences and Complexity. Vol.  24 (2), pp. 367-380.

\bibitem{van1}  van Hoeij, M. (1994). {\it  Computing Parametrizations of Rational Algebraic Curves}. In J. von zur Gathen (ed)
Proc. ISSAC 1994, pp.  187-190.


\bibitem{Vander} van der Waerden, B.L.  (1970). {\it Algebra I and II.}
Springer-Verlag, New York.

\bibitem{vdw72}  van der Waerden, B.L. (1973). {\it Einf\"{u}rung in Die Algebraischen Geometrie}. Springer Verlag, Berlin.


 \bibitem{win} Winkler, F. (1996). {\it Polynomials Algorithms in Computer
 Algebra.} Wien New York: Springer-Verlag.

\bibitem{zeng04} Zeng, Z.,  Dayton, B.H. (2004). {\it The approximate GCD of inexact polynomials part II: a multivariate algorithm}. In Proc. ISSAC 2004. pp. 320-327.

\end{thebibliography}
\end{document}